\documentclass[a4paper,10pt]{article}
\usepackage[margin=1.5in]{geometry}

\usepackage{lipsum}
\usepackage{graphicx}
\usepackage{epstopdf}
\usepackage{algorithmic}
\ifpdf
  \DeclareGraphicsExtensions{.eps,.pdf,.png,.jpg}
\else
  \DeclareGraphicsExtensions{.eps}
\fi

\usepackage{graphicx}
\usepackage{hyphenat}
\usepackage{amssymb, amsmath, amsfonts, amsopn, amsthm, array}
\usepackage{mathtools, xfrac, scalerel}
\usepackage[export]{adjustbox}
\usepackage[nocompress]{cite}

\usepackage[hidelinks]{hyperref}
\usepackage{cleveref}
\creflabelformat{equation}{#2\textup{(#1)}#3} 
\crefname{equation}{}{} 
\Crefname{equation}{}{}
\usepackage{enumitem}
\newlist{problems}{enumerate}{1}
\setlist[problems,1]{
    label={Problem~\arabic*},
    ref=\arabic*, 
    wide, labelindent=\parindent, itemsep=\parskip} 
\creflabelformat{problemsi}{#2\textup{#1}#3} 
\crefname{problemsi}{problem}{problems} 
\Crefname{problemsi}{Problem}{Problems}
\crefname{enumi}{Problem}{Problems}			

\newcommand{\opif}[2]{\underset{\scriptsize\mbox{${#2}$}}{#1}}

\DeclareMathOperator*{\diag}{diag}
\DeclareMathOperator*{\vecspan}{span}		

\newcommand{\bb}[1]{{\boldsymbol{#1}}}
\newcommand{\mc}[1]{{\mathcal{#1}}}

\usepackage{multirow}
\usepackage{eqparbox}
\newcommand{\tableIndices}[2]{{
  \begin{array}{@{}r@{}}
    \scriptstyle #2\smash{\eqmakebox[ind]{$\scriptstyle\rightarrow$}} \\[-\jot]  
    \scriptstyle #1\smash{\eqmakebox[ind]{$\scriptstyle\downarrow$}}
  \end{array}}}
\newcommand\TopSp{\rule{0pt}{2.6ex}}         
\newcommand\BotSp{\rule[-0.9ex]{0pt}{0pt}}   

\usepackage{tikz}
\usepackage{pgfplots}
\usepackage{pgfplotstable}
\usepgfplotslibrary{patchplots}
\pgfplotsset{compat=1.15}
\usepackage{arydshln}
\usetikzlibrary{math, matrix, shapes, calc, fit, backgrounds, arrows,arrows.meta, decorations.pathreplacing, decorations.markings}
\usetikzlibrary{colorbrewer}
\usepackage{listofitems,siunitx,enumerate,tabu}
\pgfdeclarelayer{background}
\pgfdeclarelayer{foreground}
\pgfsetlayers{background,main,foreground}   
\usetikzlibrary{external}
\pgfkeys{
    /tikz/external/mode=list and make
}

\newtheorem{remark}{Remark}

\title{Space-time block preconditioning for incompressible flow\date{}\thanks{Submitted to the editors 08/01/21. This work was funded by the EPSRC Centre for Doctoral Training in Industrially Focused Mathematical Modelling (EP/L015803/1), in collaboration with the Culham Centre for Fusion Energy.
BSS was supported as a Nicholas C. Metropolis Fellow under the Laboratory Directed Research and
Development program of Los Alamos National Laboratory.}}

\author{Federico Danieli\thanks{Mathematical Institute, University of Oxford, Oxford, UK 
  (\href{mailto:federico.danieli@maths.ox.ac.uk}{\tt{federico.danieli@maths.ox.ac.uk}}, \href{mailto:andy.wathen@maths.ox.ac.uk}{\tt{andy.wathen@maths.ox.ac.uk}}).}
  \and Ben S. Southworth\thanks{Los Alamos National Laboratory, Los Alamos, NM 
  (\href{mailto:southworth@lanl.gov}{\tt{southworth@lanl.gov}}), \url{http://orcid.org/0000-0002-0283-4928}}.
\and Andrew J. Wathen\footnotemark[2]}

\newenvironment{@abs}[1]{%
       \vspace{4pt}\footnotesize  \parindent 15pt {\bfseries #1. }\ignorespaces
     }
     {\par\vspace{7pt}}

\renewenvironment{abstract}{\begin{@abs}{\abstractname}}{\end{@abs}}
\newenvironment{keywords}{\begin{@abs}{\keywordsname}}{\end{@abs}}
\newenvironment{AMS}{\begin{@abs}{\AMSname}}{\end{@abs}}

\newcommand\keywordsname{Key words}
\newcommand\AMSname{AMS subject classifications}

\begin{document}

\maketitle

\begin{abstract}
  Parallel-in-time methods have become increasingly popular in the simulation
  of time-dependent numerical PDEs, allowing for the efficient use of additional
  MPI processes when spatial parallelism saturates. Most methods treat the 
  solution and parallelism in space and time separately. In contrast, all-at-once
  methods solve the full space-time system directly, largely treating time as
  simply another spatial dimension. All-at-once methods offer a number of
  benefits over separate treatment of space and time, most notably significantly
  increased parallelism and faster time-to-solution (when applicable). However,
  the development of fast, scalable all-at-once methods has largely been limited
  to time-dependent (advection-)diffusion problems. 
  This paper introduces the concept of space-time block preconditioning
  for the all-at-once solution of incompressible flow. By extending well-known
  concepts of spatial block preconditioning to the space-time setting, we
  develop a block preconditioner whose application requires the solution of a space-time
  (advection-)diffusion equation in the velocity block, coupled with a
  pressure Schur complement approximation consisting of independent spatial solves at
  each time-step, and a space-time matrix-vector multiplication. The new method is tested
  on four classical models in incompressible flow. Results indicate
  perfect scalability in refinement of spatial and temporal mesh spacing, perfect
  scalability in nonlinear Picard iterations count when applied to a nonlinear
  Navier-Stokes problem, and minimal overhead in terms of number of preconditioner
  applications compared with sequential time-stepping.
\end{abstract}




\begin{section}{Introduction}
\label{sec::intro}
{Parallel-in-Time} (PinT) algorithms have received significant attention in
recent years because of their potential to expose further concurrency in the
numerical solution of {partial differential equations} (PDEs). At some number of
MPI processes, traditional (sequential) time-stepping methods that rely on
spatial parallelism reach saturation, that is, further increasing the number of
MPI processes no longer reduces the total time to solution. Due to the
increasing number of processors available for computation, it is then desirable
to also utilise parallelism in the temporal dimension. One of the original PinT
algorithms is \emph{Parareal} \cite{pararealOriginal,pararealGander}, introduced
by Lions, Maday, and Turinici in 2001, which is essentially a two-level
multigrid method in time. Since then, a wide range of approaches for PinT have
been developed. For a complete overview of PinT algorithms, we point to the
review by Gander \cite{PinTgander}. Broadly, we will split time-parallelisation techniques
into two categories: (i) \emph{PinT} methods (excusing the slight abuse of
terminology), which treat the temporal dimension and its parallelisation
separately from the spatial dimensions, and (ii) \emph{all-at-once} methods,
which treat time and space \emph{together}, and achieve parallelisation by
solving the full monolithic system representing the space-time discretisation of
a time-dependent PDE.

A variety of all-at-once methods have been proposed \cite{ganderAAO,ganderAAOnonlin,
pestanaAAO,goddardAAO,LiuAAO}. Some of the most successful ones (particularly for the heat
equation) are rooted on multigrid techniques, and directly apply geometric multigrid (GMG)
or algebraic multigrid (AMG) to full space-time linear systems \cite{Falgout.2017,
Horton.1995,sivas2020air}. Such an approach naturally allows for
adaptive mesh refinement in both space and time \cite{sivas2020air}, and has
shown to offer significantly more parallelism and faster time-to-solution
\cite{Falgout.2017} compared with PinT methods that treat space and time separately
such as Parareal and multigrid reduction in time (MGRIT) \cite{MGRIToriginal}.
While in general multigrid approaches have been successful for the space-time parallelisation
of single-variable (advection-)diffusion equations, to our knowledge no multigrid methods
have been developed and proven effective on systems of PDEs in space and time.
Indeed, effectively handling at the same time both the complications introduced by inter-variable
coupling in systems of PDEs, as well as the purely advective term in one dimension representing a time
derivative, is beyond the capabilities of the current state-of-the-art.

Flow problems are a particularly appealing target for time-parallelisation because of
their long time dynamics and the high computational cost typically associated
with their solution. See \cite{wang2013towards} for a discussion on the need
for time parallelism in numerical fluid dynamics on emerging architectures. In
this paper, we focus on incompressible flow. Parareal
has been applied to incompressible flows in a number of works, e.g.,
\cite{miao2019convergence,trindade2006parallel,trindade2004parallel,
croce2014parallel,pararealNS}, which have generally shown successful results. 
Save for few exceptions, in most of the articles listed above the
treatment of space and time is still kept separate, and one does not
get the same benefits that can arise from a true all-at-once approach.
A number of papers have also considered the all-at-once \emph{discretisation}
of incompressible flow using space-time finite elements (for example, see
\cite{rhebergen2013space, rhebergen2012space}), but to our knowledge, none of
these papers have developed efficient, parallel methods for the all-at-once
\emph{solution} of such space-time discretisations.

In this paper, we take a new approach and develop a \emph{space-time block
preconditioner} for the all-at-once solution of incompressible flow. Applying
our preconditioner requires solving a space-time advection-diffusion
equation involving only the velocity variable, and inverting an approximate
space-time pressure Schur complement. The latter relies on the solution of a
pressure Laplacian and mass matrix at each time-step; these are independent
for each time instant, which makes
this procedure naturally parallel in time. For the former, there already exist
efficient methods for the time-parallel solution of single-variable
advection-diffusion equations, and their development remains more tractable than for
the fully-coupled space-time system. The block separation introduced by our
preconditioner allows for great freedom in the
choice of the solver for the velocity block, and it makes it straightforward to apply,
e.g., space-time multigrid as in \cite{Falgout.2017, Horton.1995,sivas2020air}.
The proposed method takes motivation from well-developed spatial block
preconditioning techniques used in sequential time-stepping or steady-state
problems \cite{andyFIT}. The development of many such methods is based on a
commuting argument involving some spatial operators, which then allows one to find
an effective approximation of the Schur complement.
Noting that formally (that is, without boundary conditions)
spatial and temporal derivatives commute, we claim that similar principles
can be generalised to the space-time setting, where an additional
time derivative is present. Indeed, in \cref{sec::precon}, our
preconditioner is derived using multiple different frameworks that can be
seen as natural extensions from the spatial block preconditioning setting.

Similar in spirit to this paper is the approach described in \cite{StollAAO}, where
the authors propose a preconditioner for a time-dependent optimal control problem, exploiting
the block structure of the system at the space-time level. However, to our knowledge, similar
concepts have not been applied before for the purpose of time-parallelisation.
Moreover, while our focus in this manuscript is incompressible flow, we believe
the underlying principle of space-time block preconditioning can be extended
to other systems of PDEs for which effective spatial block preconditioners have
been developed for the single time-step case (or for its stationary
counterpart).

The content of this article is as follows. \Cref{sec::oseen}
introduces our target problem, together with notation that will be
used throughout the paper. A derivation of the space-time block preconditioner
representing our main contribution is given in \cref{sec::precon},
together with multiple theoretical justifications. In \cref{sec::results} we
investigate the effectiveness of the proposed preconditioner, testing its
performance on four classical models in incompressible flow. Results indicate
perfect scalability in refinement of spatial and temporal mesh spacing, perfect
scalability in nonlinear Picard iteration count when applied to a nonlinear Navier-Stokes
problem, and minimal overhead in terms of number of preconditioner applications
compared with sequential time-stepping. \Cref{sec::conclusion} summarises our
findings and discusses directions for future work.

\end{section}
\begin{section}{Problem definition and discretisation}
\label{sec::oseen}
This section introduces notation that is used throughout the manuscript. Here we
describe a time-dependent version of the Oseen equations, as well as the
discrete operators arising from its finite element discretisation.

Our target system is
\begin{equation}
  \left\{\begin{array}{lr}
  \begin{array}{rcl}
    \displaystyle\frac{\partial \bb{u}}{\partial t}(\bb{x},t) + (\bb{w}(\bb{x},t)\cdot\nabla)\bb{u}(\bb{x},t)\quad\quad\\
                               - \mu \nabla^2\bb{u}(\bb{x},t) + \nabla p(\bb{x},t) &=& \bb{f}(\bb{x},t) \BotSp\\
                                                                      \nabla\cdot\bb{u}(\bb{x},t) &=& 0 \TopSp
  \end{array}  & \quad\text{in}\;\Omega\times(T_0,T],\\
  \end{array}\right.
\label{eqn::oseen}
\end{equation}
defined on a spatial domain $\Omega\subset\mathbb{R}^d$. Its unknowns are
the $d$-dimensional vector field $\bb{u}(\bb{x},t)$, and the scalar field
$p(\bb{x},t)$, representing velocity and pressure, respectively; $\mu$ is the
viscosity coefficient, which we consider constant on $\Omega$. The vector
function $\bb{f}(\bb{x},t)$ is a given forcing term, while $\bb{w}(\bb{x},t)$
represents a given $d$-dimensional advection field. The special case where
$\bb{w}(\bb{x},t)\equiv0$ gives rise to a time-dependent version of Stokes
equations:
\begin{equation}
  \left\{\begin{array}{lr}
  \begin{array}{rcl}
    \displaystyle\frac{\partial \bb{u}}{\partial t}(\bb{x},t) - \mu \nabla^2\bb{u}(\bb{x},t)+ \nabla p(\bb{x},t) &=& \bb{f}(\bb{x},t) \BotSp\\
                                                                      \nabla\cdot\bb{u}(\bb{x},t) &=& 0 \TopSp
  \end{array}  & \quad\text{in}\;\Omega\times(T_0,T].\\
  \end{array}\right.
\label{eqn::stokes}
\end{equation}
While most of the analysis conducted in this paper considers the more general problem
\cref{eqn::oseen}, we will at times refer to \cref{eqn::stokes} for some
simplifications. Conversely, by introducing the nonlinearity $\bb{w}(\bb{x},t)\equiv\bb{u}(\bb{x},t)$
we produce the Navier-Stokes equations for incompressible flow
\begin{equation}
  \left\{\begin{array}{lr}
  \begin{array}{rcl}
    \displaystyle\frac{\partial \bb{u}}{\partial t}(\bb{x},t) + (\bb{u}(\bb{x},t)\cdot\nabla)\bb{u}(\bb{x},t)\quad\quad\\
                               - \mu \nabla^2\bb{u}(\bb{x},t) + \nabla p(\bb{x},t) &=& \bb{f}(\bb{x},t) \BotSp\\
                                                                      \nabla\cdot\bb{u}(\bb{x},t) &=& 0 \TopSp
  \end{array}  & \quad\text{in}\;\Omega\times(T_0,T].\\
  \end{array}\right.
\label{eqn::navierStokes}
\end{equation}
The Oseen system can then be interpreted as a linearisation of Navier-Stokes,
implying that the analysis conducted here can be extended to more general frameworks:
this is addressed more in detail in \cref{sec::results::NS}.

In order to provide closure to the definition of the initial-boundary value
problem, we also equip \cref{eqn::oseen} with an {initial condition}
\begin{equation}
  \bb{u}(\bb{x},T_0) = \bar{\bb{u}}^0(\bb{x})\quad\text{on }\Omega,
  \label{eqn::oseenIC}
\end{equation}
and appropriate {boundary conditions} (BC)
\begin{equation}
\begin{array}{rcll}
  \bb{u}(\bb{x},t)                                                        &=& \bb{u}_D(\bb{s},t)\quad & \text{on }\Gamma_D\times(T_0,T],\quad\text{and} \BotSp\\
  \mu\nabla\bb{u}(\bb{s},t)\,\bb{n}(\bb{s}) - p(\bb{s},t)\,\bb{n}(\bb{s}) &=& \bb{g}(\bb{s},t)  \quad & \text{on }\Gamma_N\times(T_0,T], \TopSp
\end{array}
\label{eqn::oseenBC}
\end{equation}
where we split the boundary $\partial\Omega$ into its Dirichlet and Neumann
parts, respectively, with $\partial\Omega=\Gamma_D\cup\Gamma_N$ and
$\Gamma_D\cap\Gamma_N=\varnothing$. In \cref{eqn::oseenBC},
$\bb{n}(\bb{s})\in\mathbb{R}^{d}$ represents the outward-pointing normal to the
boundary, while $\bb{u}_D(\bb{s},t)$, $\bb{g}(\bb{s},t)$, and
$\bar{\bb{u}}^0(\bb{x})$ are some known functions.

At any given instant $t\in\left[T_0,T\right]$, the velocity unknown lives in the
functional space \cite[Chap.~2.4.2]{quarterNM}
\begin{equation}
  \bb{u}(\cdot,t)\in V(\Omega) \equiv \left(H^{1}_{\Gamma_D}(\Omega)\right)^d,
  \label{eqn::funcSpaceV}
\end{equation}
including the Dirichlet condition in the form of $\bb{u}_D$, while for the pressure we have
\begin{equation}
  p(\cdot,t)\in Q(\Omega)\equiv L^{2}(\Omega),\quad\text{or}\quad Q(\Omega)\equiv L^2_0(\Omega)=\left\{q(\bb{x})\in L^2(\Omega):\int_\Omega q(\bb{x})\,d\bb{x}=0\right\},
  \label{eqn::funcSpaceP}
\end{equation}
where the zero-mean valued counterpart, $L^2_0(\Omega)$, is chosen if Dirichlet
boundary conditions are imposed everywhere on the boundary of the domain,
$\Gamma_D\equiv\partial\Omega$. In either case, having introduced
\cref{eqn::funcSpaceV,eqn::funcSpaceP}, we can provide the weak formulation
\cite[Chap.~17.2]{quarterNM} of \cref{eqn::oseen}: we seek $\bb{u}(\bb{x},t)$
and $p(\bb{x},t)$ such that
\begin{equation}\left\{
  \begin{array}{rcll}
    \displaystyle\frac{\partial}{\partial t}\left\langle\bb{u},\bb{v}\right\rangle+\mu \left\langle\nabla\bb{u},\nabla\bb{v}\right\rangle\quad\quad\quad & & &\\
    \displaystyle+\left\langle(\bb{w}\cdot\nabla)\bb{u},\bb{v}\right\rangle-\left\langle p,\nabla\cdot\bb{v}\right\rangle&=&\left\langle\bb{f},\bb{v}\right\rangle
     + \displaystyle\int_{\Gamma_N}\bb{g}\cdot\bb{v} &\quad\forall\bb{v}(\bb{x})\in V(\Omega)\BotSp\\
    \displaystyle- \left\langle q,\nabla\cdot\bb{u}\right\rangle &=& 0 &\quad\forall q(\bb{x})\in Q(\Omega)\TopSp\\
  \end{array},\right.
\label{eqn::oseenWeak}
\end{equation}
for each instant $t\in(T_0,T]$, where $\left\langle
\bb{a}(\bb{x}),\bb{b}(\bb{x})\right\rangle =
\int_{\Omega}\bb{a}(\bb{x})\cdot\bb{b}(\bb{x})\,d\bb{x}$ denotes the spatial scalar
product (notice we dropped the variables' dependence on $\bb{x}$ and $t$ in order
to contract notation).

\paragraph{Space-time discretisation}
The temporal derivative in \cref{eqn::oseenWeak} is discretised using {implicit Euler}.
To this end, we introduce a uniformly-spaced grid, $\{t_k\}$, over the temporal domain, where
\begin{equation}
  t_k = T_0 + k\Delta t,\qquad\text{with}\qquad\Delta t = (T-T_0)/N_t \qquad\text{and}\qquad k=0,\dots,N_t,
  \label{eqn::timeDisc}
\end{equation}
at whose nodes we seek to recover our numerical solution. Following a
{finite element} approach \cite[Chap.~4]{quarterNM} \cite[Chap.~10]{andyFIT}, we discretise
\cref{eqn::oseenWeak} in space by introducing a triangularisation of $\Omega$ of characteristic
length $\Delta x$, on which we build finite-dimensional approximations to
both $V(\Omega)$ and $Q(\Omega)$, named $V_h(\Omega)$ and $Q_h(\Omega)$.
These are identified by their basis functions
\begin{subequations}
\begin{align}
  V_h(\Omega)&\coloneqq\vecspan\{\bb{\phi}_0(\bb{x}),\dots,\bb{\phi}_{N_{\bb{u}}-1}(\bb{x})\},\quad\text{and}\quad\label{eqn::VQbases::V}\\
  Q_h(\Omega)&\coloneqq\vecspan\{\psi_0(\bb{x}),\dots,\psi_{N_p-1}(\bb{x})\}\label{eqn::VQbases::Q},
\end{align}
\label{eqn::VQbases}%
\end{subequations}
with $N_{\bb{u}}$ and $N_p$ being the number of degrees of freedom associated
with the velocity and pressure variables, respectively, so that vectors
$\bb{v}\in\mathbb{R}^{N_{\bb{u}}}$ and $\bb{q}\in\mathbb{R}^{N_{p}}$ correspond
to functions in $V_h(\Omega)$ and $Q_h(\Omega)$ via
\begin{equation}
  \bb{v}(\bb{x}) =\sum_{j=0}^{N_{\bb{v}}-1}\left[\bb{v}\right]_j\bb{\phi}_j(\bb{x}) \qquad\text{and}\qquad q(\bb{x}) = \sum_{j=0}^{N_{p}-1}\left[\bb{q}\right]_j\psi_j(\bb{x}),
  \label{eqn::FEvec2fun}
\end{equation}
respectively. In the remainder of the manuscript, we will refer somewhat
interchangeably to functions in our discrete spaces and their vector
representations according to the bases \eqref{eqn::VQbases}. With this, we can
assemble the {Galerkin matrices}, namely the discrete counterparts of the
bilinear forms appearing in \eqref{eqn::oseenWeak}. In particular, we define the
{mass} and {stiffness} matrices for the velocity variable, respectively,
\begin{subequations}
  \begin{align}
    \left[\mc{M}_{\bb{u}}\right]_{m,n=0}^{N_{\bb{u}}-1}&\coloneqq\int_{\Omega}      \bb{\phi}_{m}(\bb{x})\cdot       \bb{\phi}_{n}(\bb{x})\,d\bb{x}\quad\text{and}\quad\label{eqn::massStiffV::mass}\\
    \left[\mc{A}_{\bb{u}}\right]_{m,n=0}^{N_{\bb{u}}-1}&\coloneqq\int_{\Omega}\nabla\bb{\phi}_{m}(\bb{x})\colon\nabla\bb{\phi}_{n}(\bb{x})\,d\bb{x}\label{eqn::massStiffV::stiff},
  \end{align}
\label{eqn::massStiffV}%
\end{subequations}
the discrete (negative) {divergence} matrix, responsible for the
coupling between velocity and pressure,
\begin{equation}
  \left[\mc{B}\right]_{m,n=0}^{N_{p}-1,N_{\bb{u}}-1} \coloneqq -\int_{\Omega}{\psi}_{m}(\bb{x})\nabla\cdot \bb{\phi}_{n}(\bb{x})\,d\bb{x},
\end{equation}
and discretisations of the advection operator for each temporal instant $t_k$ in \cref{eqn::timeDisc}:
\begin{equation}
  \left[\mc{W}_{\bb{u},k}\right]_{m,n=0}^{N_{\bb{u}}-1} \coloneqq \int_{\Omega} \left(\left(\bb{w}(\bb{x},t_k)\cdot\nabla\right)\bb{\phi}_{n}(\bb{x})\right)\cdot\bb{\phi}_{m}(\bb{x})\,d\bb{x}.
\label{eqn::advV}
\end{equation}

Combining spatial and temporal discretisations gives rise to the following discrete
linear system approximating \cref{eqn::oseenWeak}:
\begin{equation}
  \left\{\begin{array}{rl}
    \displaystyle\frac{1}{\Delta t}\mc{M}_{\bb{u}}\left(\bb{u}^{k}-\bb{u}^{k-1}\right)+\mu\mc{A}_{\bb{u}}\bb{u}^{k}+\mc{W}_{\bb{u},k}\bb{u}^{k}+\mc{B}^T\bb{p}^{k} &= \bb{f}^{k}\\
                                                                                                                                                  \mc{B}\bb{u}^{k} &= \bb{0}
  \end{array}\right.,\quad k=1,\dots N_t,
  \label{eqn::oseenSingleStep}
\end{equation}
where $\bb{0}$ is an all-zero vector of size $N_p$. Here,
$\left[\bb{u}^k\right]_j$ and $\left[\bb{p}^k\right]_j$ represent the
coefficients corresponding to the $j$-th basis function for the velocity and
pressure variables in \cref{eqn::VQbases::V,eqn::VQbases::Q}, respectively,
evaluated at the $k$-th temporal instant; analogously, $[\bb{f}^k]_j$ identifies
the $L^2$-projection of the right-hand side of \cref{eqn::oseenWeak} on the
$j$-th basis function for $V_h(\Omega)$, at the instant $t_k$:
\begin{equation}
  [\bb{f}^k]_j \coloneqq \int_{\Omega}\bb{f}(\bb{x},t_k)\cdot\bb{\phi}_j(\bb{x})\,d\bb{x} + \int_{\Gamma_N}\bb{g}(\bb{s},t_k)\cdot\bb{\phi}_j(\bb{s})\,d\bb{s}.
  \label{sec::oseenRhs}
\end{equation}

We now have all the ingredients necessary to assemble the monolithic space-time system corresponding to \eqref{eqn::oseenSingleStep}.
To simplify notation, we introduce the operators
\begin{equation}
  \mc{F}_{\bb{u},k} \coloneqq \frac{\mc{M}_{\bb{u}}}{\Delta t} + \mathcal{W}_{\bb{u},k}  + \mu\mc{A}_{\bb{u}},\quad\quad k=1,\dots N_t.
  \label{eqn::VRAD}
\end{equation}
We can then arrange all the equations of \cref{eqn::oseenSingleStep} in a single block matrix, as
\begin{equation}
\begin{tikzpicture}[
      baseline=(current  bounding  box.center), 
      Highlight/.style={
          draw,
          densely dotted,
          rounded corners=2pt,
      },
  ]
  \matrix[matrix of nodes,align=center,nodes in empty cells,left delimiter={[},right delimiter ={]}] at (0,0) (A){ 
  $\mc{F}_{\bb{u},1}               $&$\mc{B}^T$&          &     &                     &          \\
  $\mc{B}                          $&          &          &     &                     &          \\
  $\frac{\mc{M}_{\bb{u}}}{\Delta t}$&          &$\;\ddots$&     &                     &          \\
                                    &          &          &\quad&                     &          \\
                                    &          &$\;\ddots$&     &$\mc{F}_{\bb{u},N_t}$&$\mc{B}^T$\\
                                    &          &          &     &$\mc{B}             $&          \\
  };
  \matrix[matrix of nodes,align=center,nodes in empty cells,left delimiter={[},right delimiter ={]},inner xsep=-1pt,anchor=west,right=17pt] at (A.east) (X){ 
  $\bb{u}^1    $\\
  $\bb{p}^1    $\\
  $\vdots      $\\[5pt]
  $            $\\
  $\bb{u}^{N_t}$\\
  $\bb{p}^{N_t}$\\
  };
  \node[anchor=west,right=7pt] at (X.east) (E) {$=$};
  \matrix[matrix of nodes,align=center,nodes in empty cells,left delimiter={[},right delimiter ={]},inner xsep=-1pt,anchor=west,right=7pt] at (E.east) (B){ 
  $\bb{f}^1    $\\
  $\bb{0}      $\\
  $\vdots      $\\[5pt]
  $            $\\
  $\bb{f}^{N_t}$\\
  $\bb{0}      $\\
  };
  \node[anchor=west,right=7pt] at (B.east) (D) {$,$};
  \draw[Highlight](X-2-1.south -| X-5-1.east) -- (X-2-1.south -| X-5-1.west);
  \draw[Highlight](X-5-1.north east) -- (X-5-1.north west);
  \draw[Highlight](B-2-1.south -| B-5-1.east) -- (B-2-1.south -| B-5-1.west);
  \draw[Highlight](B-5-1.north east) -- (B-5-1.north west);
  \draw[Highlight](A-1-2.north east) rectangle (A-2-1.south -| A-1-1.west);
  \draw[Highlight](A-3-1.north -| A-1-1.west) rectangle (A-4-1.south -| A-1-2.east);
  \draw[Highlight](A-5-6.north east) rectangle (A-6-5.south -| A-5-5.west);
\end{tikzpicture}
\label{eqn::STStepDisc}
\end{equation}
where, with a slight abuse of notation, we re-define $\bb{f}^1$ from \cref{sec::oseenRhs} in order to include  the initial condition \cref{eqn::oseenIC},
\begin{equation}
  \left[\bb{f}^1\right]_j \coloneqq \int_{\Omega}\left(\bb{f}(\bb{x},t_1)+\frac{1}{\Delta t}\bar{\bb{u}}^0(\bb{x})\right)\cdot\bb{\phi}_j(\bb{x})\,d\bb{x} + \int_{\Gamma_N}\bb{g}(\bb{s},t_k)\cdot\bb{\phi}_j(\bb{s})\,d\bb{s}.
\end{equation}

Notice that the system in \cref{eqn::STStepDisc} is of the block
lower-triangular structure typical of space-time discretisations
\cite{pararealGander,MGRIToriginal}; moreover, since we use a \emph{one}-step
method to approximate the time derivative, the system is made of
\emph{two} block diagonals, with the $2\times2$ blocks composing them
highlighted in \cref{eqn::STStepDisc}. Usually, the solution to
\cref{eqn::oseenSingleStep} is recovered via a time-stepping routine, which
corresponds to solving \cref{eqn::STStepDisc} directly via (block) {forward
substitution}. There are several methods to accelerate this
procedure, mostly focusing on accelerating the inversion of the main diagonal
blocks
\begin{equation}
  \mc{A}_k\coloneqq\left[\begin{array}{cc}
    \mc{F}_{\bb{u},k} & \mc{B}^T\\
    \mc{B}
  \end{array}\right],\qquad k=1,\dots,N_t,
  \label{eqn::oseenBlock}
\end{equation}
by designing an optimal preconditioner (see for example
\cite[Chap.~17.8]{quarterNM}, or \cite[Chap.~4]{andyFIT} for an overview). The \emph{pressure convection-diffusion} (PCD) preconditioner is
one of the most successful among these \cite{PCDoriginal}, preconditioning the
block structure of \cref{eqn::oseenBlock} by finding an accurate and computable
approximation to the {pressure Schur complement}. Our work stems from the
same principles as PCD, but considers the whole space-time system
\cref{eqn::STStepDisc}, rather than just the spatial diagonal blocks
\cref{eqn::oseenBlock}. Motivation and derivation of our preconditioner are
described in the following section.

\end{section}
\begin{section}{Space-time block-preconditioning}
\label{sec::precon}
Rather than inverting \cref{eqn::STStepDisc} using sequential time-stepping
(i.e., forward substitution), here we tackle the solution of the full space-time
system \cref{eqn::STStepDisc} all at once by using preconditioned GMRES
\cite{GMRESoriginal}, with an appropriate space-time block preconditioner. In
particular, we focus on developing a preconditioner that treats space and time
together and in parallel, which may yield faster time-to-solution than block
forward substitution when sufficient processors are available. In the remainder
of this section, we show how we can effectively build on concepts developed in
the study of the single time-step matrix \cref{eqn::oseenBlock} to design
preconditioners that are effective in the full space-time setting.

To develop our preconditioner, we start by reordering the space-time
operator in \cref{eqn::STStepDisc} to take the form
\begin{equation}
\begin{tikzpicture}[
      baseline=(current  bounding  box.center), 
      Highlight/.style={
          draw,
      },
  ]
  \matrix[matrix of nodes,align=center,nodes in empty cells,left delimiter={[},right delimiter ={]}] at (0,0) (A){ 
  $\mc{F}_{\bb{u},1}                $&$          $&$                   $&$\mc{B}^T$&$          $&$        $\\[-7pt]
  $-\frac{\mc{M}_{\bb{u}}}{\Delta t}$&$\;\ddots\;$&$                   $&$        $&$\;\ddots\;$&$        $\\[-7pt]
  $                                 $&$\;\ddots\;$&$\mc{F}_{\bb{u},N_t}$&$        $&$          $&$\mc{B}^T$\\
  $\mc{B}                           $&$          $&$                   $&$        $&$          $&$        $\\
  $                                 $&$\;\ddots\;$&$                   $&$        $&$          $&$        $\\
  $                                 $&$          $&$\mc{B}             $&$        $&$          $&$        $\\
  };
  \matrix[matrix of nodes,align=center,nodes in empty cells,left delimiter={[},right delimiter ={]},inner xsep=-1pt,inner ysep=2.5pt,anchor=west,right=17pt] at (A.east) (X){ 
  $\bb{u}^1    $\\
  $\vdots      $\\[4pt]
  $\bb{u}^{N_t}$\\
  $\bb{p}^1    $\\
  $\vdots      $\\
  $\bb{p}^{N_t}$\\
  };
  \node[anchor=west,right=7pt] at (X.east) (E) {$=$};
  \matrix[matrix of nodes,align=center,nodes in empty cells,left delimiter={[},right delimiter ={]},inner xsep=-1pt,inner ysep=3pt,anchor=west,right=7pt] at (E.east) (B){ 
  $\bb{f}^1    $\\[-4pt]
  $\vdots      $\\
  $\bb{f}^{N_t}$\\[3pt]
  $\bb{0}      $\\
  $\vdots      $\\
  $\bb{0}      $\\
  };
  \node[anchor=west,right=7pt] at (B.east) (D) {$.$};
  \draw[Highlight](X-3-1.south -| X-3-1.east) -- (X-3-1.south -| X-3-1.west);
  \draw[Highlight](B-3-1.south -| B-3-1.east) -- (B-3-1.south -| B-3-1.west);
  \draw[Highlight]($ (A.north) + (10pt,0) $) -- ($ (A.south) + (10pt,0) $);
  \draw[Highlight](A.east)  -- (A.west);
\end{tikzpicture}
\label{eqn::STdisc}
\end{equation}
This rearrangement favours a separation of variables based on the physical
unknown they refer to, rather than on the individual time-steps. We further
denote the blocks composing \cref{eqn::STdisc} by
\begin{equation}
  F_{\bb{u}}\coloneqq\left[\begin{array}{ccc}
    \mc{F}_{\bb{u},1}                &      &                   \\
    -\frac{\mc{M}_{\bb{u}}}{\Delta t}&\ddots&                   \\
                                     &\ddots&\mc{F}_{\bb{u},N_t}
  \end{array}\right],\qquad\text{and}\qquad B \coloneqq \diag_{N_t}(\mc{B}),
  \label{eqn::STFu} 
\end{equation}
where $\diag_N(*)$ denotes a block diagonal matrix with $N$ blocks containing $(*)$.
As in the single time-step case \cite[Chap.~9.2]{andyFIT}, we exploit the block structure of \cref{eqn::STdisc}, and look for a right preconditioner for GMRES in the following \emph{block upper triangular} form:
\begin{equation}
  P_T      \coloneqq \left[\begin{array}{cc}F_\bb{u}      & B^T                   \\  & -X      \end{array}\right] \quad\Longleftrightarrow\quad
  P_T^{-1} \coloneqq \left[\begin{array}{cc}F_\bb{u}^{-1} & F_\bb{u}^{-1}B^TX^{-1}\\  & -X^{-1} \end{array}\right],
  \label{eqn::PT}
\end{equation}
where $X$ approximates the space-time pressure Schur complement
\begin{equation}
  X \approx S_p\coloneqq BF_\bb{u}^{-1}B^T.
  \label{eqn::schur}
\end{equation}
Convergence of GMRES applied to the full space-time system and preconditioned
by $P_T^{-1}$ is then exactly defined by convergence of GMRES applied to the
Schur-complement problem, preconditioned by $X^{-1}$ \cite{southworth2020FP}.

The key advantage of the reordering in \cref{eqn::STdisc} and the block 
preconditioning discussed above is that we have effectively decoupled the imposition
of the Lagrangian constraint from the solution of the PDE. Applying the preconditioner
then requires inverting $F_\bb{u}$ (that is, time-integrating the velocity
field), and applying an approximate Schur complement inverse \cref{eqn::schur}.
The former represents the discretisation of a single-variable parabolic PDE,
which is much simpler to solve in parallel in space-time than the fully-coupled
system \cref{eqn::STdisc}; indeed parallel-in-time and -space-time algorithms have shown to behave well
on such problems \cite{Falgout.2017,sivas2020air}. Moreover, we claim that the pressure Schur
complement \cref{eqn::schur} can naturally be approximated in a time-parallel
fashion. To derive
such an approximation, we follow two approaches which have been inspirational in
designing efficient preconditioners for the single time-step case.

\begin{remark}
One could also consider a block diagonal preconditioner, but due to the
nonsymmetry of \cref{eqn::STFu,eqn::schur}, MINRES \cite{MINRES} cannot be
applied, thus negating one of the primary benefits of using a preconditioner in this form.
With GMRES, block diagonal preconditioners are typically expected to require
twice as many iterations \cite{southworth2020FP} as block triangular ones, and in
this case offer only a marginal reduction in computational cost. Block lower
triangular preconditioning is another option, but numerical tests (not
included here) have favoured using an upper triangular preconditioner.
The latter is also more suited for right-preconditioning
\cite{southworth2020FP}, which is required by flexible GMRES \cite{FGMRES} (and
flexible Krylov methods are necessary if we use an inner Krylov method to
approximately invert $F_\bb{u}$). For these reasons, $P_T$ \cref{eqn::PT}
remains our preconditioner of choice.
\end{remark}

\subsection{Small commutator approach}
\label{sec::precon::smallComm}
In \cite{PCDoriginal}, an effective preconditioning strategy for \cref{eqn::oseenBlock} is recovered by
starting from an assumption on the commuting of certain operators. Although our
Schur complement \cref{eqn::schur} differs from the case considered there, as it
involves a time-derivative, we can still follow a similar line of reasoning. Let us
begin by characterising our space-time Schur complement in more detail. Firstly, we can
explicitly write $F_\bb{u}^{-1}$ in the following form:\footnote{If we were
solving the time-dependent Stokes equations (that is, if
$\bb{w}(\bb{x},t)\equiv0$), the structure of \cref{eqn::BiDiagInverse} would be
further simplified, as $F_\bb{u}$ would be {block Toeplitz}; the presence
of a time-dependent advection field, instead, causes the main diagonal block to
vary down the diagonal. A similar disruption would be caused by using a
non-uniform temporal grid when discretising Stokes: in that context, the
difference between the blocks down the same diagonal is given by the factor
$\Delta t_k$ by which we divide the mass matrix. A similar analysis as the one
reported here can be applied to that case as well.}
\begin{equation}
F_\bb{u}^{-1} = \left[\begin{array}{ccc}
  \mc{F}_{\bb{u},1}^{-1}                                                                   &                      &\\
  \mc{F}_{\bb{u},2}^{-1}\left(\frac{\mc{M}_{\bb{u}}}{\Delta t}\mc{F}_{\bb{u},1}^{-1}\right)&\mc{F}_{\bb{u},2}^{-1}&\\
  \mc{F}_{\bb{u},3}^{-1}\left(\frac{\mc{M}_{\bb{u}}}{\Delta t}\mc{F}_{\bb{u},2}^{-1}
                              \frac{\mc{M}_{\bb{u}}}{\Delta t}\mc{F}_{\bb{u},1}^{-1}\right)&\mc{F}_{\bb{u},3}^{-1}\left(\frac{\mc{M}_{\bb{u}}}{\Delta t}\mc{F}_{\bb{u},2}^{-1}\right) & \ddots\\
  \vdots &\ddots&\ddots\\
\end{array}\right].
\label{eqn::BiDiagInverse}
\end{equation}
Multiplying on the left and right by the block diagonal matrices $B$ and $B^T$,
respectively, results in inner operators (i.e., blocks) of the form
\begin{equation}
  \mc{B}\mc{F}_{\bb{u},k}^{-1}\left(\prod_{i=1}^{j}\left(\frac{\mc{M}_{\bb{u}}}{\Delta t}\mc{F}_{\bb{u},k-i}^{-1}\right)\right)\mc{B}^T,
  \quad\text{with}\quad\left\{\begin{array}{l} k=1,\dots,N_t\\j=0,\dots,k-1\end{array}\right.,
  \label{eqn::FinvBlock}
\end{equation}
where the product is assumed to expand towards the right, i.e.,
$\prod_{i=0}^{k}a_i = a_{0}\prod_{i=1}^{k}a_i$, and returns identity if the
subscript is larger than the superscript.
We seek to adequately approximate the operators appearing in
\cref{eqn::FinvBlock}: to do so, we follow \cite{PCDoriginal} as well as
\cite[Chap.~9.2]{andyFIT}, and make the assumption that the
reaction-advection-diffusion operator and the gradient operator approximately
commute, that is:
\begin{equation}
  \nabla\left(\frac{\mc{I}}{\Delta t}+\left(\bb{w}(\bb{x},t)\cdot\nabla\right) - \mu\nabla^2\right)
       -\left(\frac{\mc{I}}{\Delta t}+\left(\bb{w}(\bb{x},t)\cdot\nabla\right) - \mu\nabla^2\right)\nabla\approx0.
  \label{eqn::commutatorStrong}
\end{equation}
Notice that the ``reaction'' term stems from the presence of the temporal derivative,
and is mimicked by the scaled identity operator $\mc{I}/\Delta t$, which automatically
commutes with the gradient operator. To provide a discrete analogous of
\cref{eqn::commutatorStrong}, let us define, for each $\mc{F}_{\bb{u},k}$ in \cref{eqn::VRAD},
an equivalent reaction-advection-diffusion operator acting on the pressure field:
\begin{equation}
\begin{split}
  \mc{F}_{p,k}\coloneqq&\frac{\mc{M}_{p}}{\Delta t} + \mathcal{W}_{p,k}  + \mu\mc{A}_{p}, \qquad k = 1,\dots,N_t,
\end{split}
\label{eqn::PRAD}
\end{equation}
where $\mc{M}_p$, $\mc{A}_p$, and $\mc{W}_{p,k}$ play the same role as the mass,
stiffness, and advection matrices
(\cref{eqn::massStiffV::mass}, \cref{eqn::massStiffV::stiff}, and \cref{eqn::advV}), but for the
pressure variable, namely
\begin{subequations}
  \begin{align}
    \left[\mc{M}_{p  }\right]_{m,n=0}^{N_{p}-1} &\coloneqq \int_{\Omega}       \psi_{m}(\bb{x})           \psi_{n}(\bb{x})\,d\bb{x},\label{eqn::massStiffP::mass}\\
    \left[\mc{A}_{p  }\right]_{m,n=0}^{N_{p}-1} &\coloneqq \int_{\Omega} \nabla\psi_{m}(\bb{x})\cdot\nabla\psi_{n}(\bb{x})\,d\bb{x},\quad\text{and}\label{eqn::massStiffP::stiff}\\
    \left[\mc{W}_{p,k}\right]_{m,n=0}^{N_{p}-1} &\coloneqq \int_{\Omega} \left(\bb{w}(\bb{x},t_k)\cdot\nabla\right)\psi_{n}(\bb{x})\psi_{m}(\bb{x})\,d\bb{x},\quad\text{for}\quad k = 1,\dots,N_t\label{eqn::massStiffP::adv}.
  \end{align}
  \label{eqn::massStiffP}%
\end{subequations}
More details on how to assemble these operators are provided in \cref{sec::precon::Apdef}.
Translating the commutation assumption to the discrete form introduced in
\cref{sec::oseen}, we have
\begin{equation}
  \left(\mc{M}_{\bb{u}}^{-1}\mc{B}^T         \right)\left(\mc{M}_{p     }^{-1}\mc{F}_{p,j}\right) - 
  \left(\mc{M}_{\bb{u}}^{-1}\mc{F}_{\bb{u},j}\right)\left(\mc{M}_{\bb{u}}^{-1}\mc{B}^T    \right)\approx 0
\label{eqn::commutator}
\end{equation}
for each $ j=1,\dots,N_t$. Here, re-scaling by the corresponding mass matrices
is necessary, since the Galerkin operators have the effect of mapping functions
to \emph{functionals}, integrating over the whole spatial domain.

To get closer to our desired operators \cref{eqn::FinvBlock}, we left-multiply
\cref{eqn::commutator} by $\mc{F}_{\bb{u},j}^{-1}\mc{M}_{\bb{u}}$ and
right-multiply by $\mc{F}_{p,j}^{-1}\mc{M}_{p}$, thus obtaining
\begin{equation}
  \mc{F}_{\bb{u},j}^{-1}\mc{B}^T\approx
  \mc{M}_{\bb{u}  }^{-1}\mc{B}^T\mc{F}_{p,j}^{-1}\mc{M}_{p}.
  \label{eqn::commutatorfirstblock}
\end{equation}
Continuing to left-multiply by
$\mc{F}_{\bb{u},j+i}^{-1}\frac{\mc{M}_{\bb{u}}}{\Delta t}$, increasing $i$ until
we reach $j+i=k$, using \cref{eqn::commutatorfirstblock} for each $j+i$, and
finally left-multiplying by $\mc{B}$, gives us
\begin{equation}
\begin{split}
  & \mc{F}_{\bb{u},j+1}^{-1}\left(\frac{\mc{M}_{\bb{u}}}{\Delta t}\mc{F}_{\bb{u},j}^{-1}\right)\mc{B}^T\approx
  \underbrace{\frac{\mc{F}_{\bb{u},j+1}^{-1}}{\Delta t}\mc{B}^T}_{\approx\mc{M}_{\bb{u}}^{-1}\mc{B}^T\mc{F}_{p,j+1}^{-1}\frac{\mc{M}_{p}}{\Delta t}}\mc{F}_{p,j}^{-1}\mc{M}_{p}\\
  \Longrightarrow\quad& \mc{B}\mc{F}_{\bb{u},k}^{-1}\left(\prod_{i=j}^{k-1}\left(\frac{\mc{M}_{\bb{u}}}{\Delta t}\mc{F}_{\bb{u},i}^{-1}\right)\right)\mc{B}^T\approx
  \mc{B}\mc{M}_{\bb{u}}^{-1}\mc{B}^T\mc{F}_{p,k}^{-1}\left(\prod_{i=j}^{k-1}\left(\frac{\mc{M}_{p}}{\Delta t}\mc{F}_{p,i}^{-1}\right)\right)\mc{M}_{p},
\end{split}
\label{eqn::commutatorallblocks}
\end{equation}
for any desired index $k=1,\dots,N_t$. This produces an approximation for each of the
operators in \cref{eqn::FinvBlock}, and consequently identifies the following
candidate for a reasonable approximation of the whole space-time pressure Schur
complement \cref{eqn::schur}:
\begin{equation}
  X \coloneqq \diag_{N_t}(\mc{B}\mc{M}_{\bb{u}}^{-1}\mc{B}^T)\left[\begin{array}{ccc}
  \mc{F}_{p,1}^{-1}                                                         &                 &\\
  \mc{F}_{p,2}^{-1}\left(\frac{\mc{M}_{p}}{\Delta t}\mc{F}_{p,1}^{-1}\right)&\mc{F}_{p,2}^{-1}&\\
  \mc{F}_{p,3}^{-1}\left(\frac{\mc{M}_{p}}{\Delta t}\mc{F}_{p,2}^{-1}
                         \frac{\mc{M}_{p}}{\Delta t}\mc{F}_{p,1}^{-1}\right)&\mc{F}_{p,3}^{-1}\left(\frac{\mc{M}_{p}}{\Delta t}\mc{F}_{p,2}^{-1}\right) & \ddots\\
  \vdots &\ddots&\ddots\\
\end{array}\right]\diag_{N_t}(\mc{M}_p).
  \label{eqn::XXfirstDef}
\end{equation}
The expression above can be simplified by noticing that the central matrix
shares the same structure as \cref{eqn::BiDiagInverse}. In fact, it
is exactly the inverse of the block bi-diagonal matrix
\begin{equation}
  F_{p}\coloneqq\left[\begin{array}{ccc}
    \mc{F}_{p,1}                &      &              \\
    -\frac{\mc{M}_{p}}{\Delta t}&\ddots&              \\
                                &\ddots&\mc{F}_{p,N_t}
  \end{array}\right],
  \label{eqn::STFp}
\end{equation}
which can be interpreted as a system stemming from the implicit Euler discretisation
of a time-dependent PDE, just like its velocity counterpart \cref{eqn::STFu}.
Moreover, by making use of the discrete $\inf$-$\sup$ condition, it has been
shown \cite[Chap.~3.5]{andyFIT} that the operator
$\mc{B}\mc{M}_{\bb{u}}^{-1}\mc{B}^T$ in \cref{eqn::XXfirstDef} is spectrally
equivalent to a discrete Laplacian operator acting on the pressure field. In
principle, this operator can differ from the one appearing in \cref{eqn::PRAD},
hence we denote it as
\begin{equation}
  \tilde{\mc{A}}_p \approx \mc{B}\mc{M}_{\bb{u}}^{-1}\mc{B}^T,
  \label{eqn::Plaplacian}
\end{equation}
and delay a more in-depth discussion on its definition to
\cref{sec::precon::Apdef}. Furthermore, to simplify notation, we define the
following block diagonal matrices
\begin{equation}
  \quad M_p\coloneqq\diag_{N_t}(\mc{M}_p)\qquad\text{and}\qquad A_p\coloneqq\diag_{N_t}\left(\tilde{\mc{A}}_p\right),
  \label{eqn::MpdiagApdiag}
\end{equation}
which ultimately allows us to write our approximation to the space-time
pressure Schur complement \cref{eqn::XXfirstDef} as
\begin{equation}
  X \coloneqq A_p F_p^{-1}M_p \quad\Longleftrightarrow\quad X^{-1} \coloneqq M_p^{-1} F_p A_p^{-1}.
  \label{eqn::XXdef}
\end{equation}

We point out that the derivation in this section has been outlined with the
purpose of drawing a connection with the commuting assumption typically made in
the design of Schur complement approximations for the steady-state or time-stepping
frameworks. An analogous approximation can be derived by directly considering the
commutator between the gradient and the whole parabolic operator,
\begin{equation}
  \nabla\left(\frac{\partial}{\partial t}+\left(\bb{w}(\bb{x},t)\cdot\nabla\right) - \mu\nabla^2\right)
       -\left(\frac{\partial}{\partial t}+\left(\bb{w}(\bb{x},t)\cdot\nabla\right) - \mu\nabla^2\right)\nabla\approx0,
  \label{eqn::commutatorStrongST}
\end{equation}
(which is equivalent to \cref{eqn::commutatorStrong}, given how temporal and spatial
derivatives automatically commute). In fact, its discrete counterpart is given by
\begin{equation}
  (M_{\bb{u}}^{-1}B^T)(M_p^{-1}F_p) - (M_{\bb{u}}^{-1}F_{\bb{u}})(M_{\bb{u}}^{-1}B^T)\approx0,
  \label{eqn::commutatorST}
\end{equation}
with $M_{\bb{u}}\coloneqq\diag_{N_t}(\mc{M}_{\bb{u}})$. Left-multiplying by
$BF_{\bb{u}}^{-1}M_{\bb{u}}$ and right-multiplying by $F_p^{-1}M_p$, we recover
\begin{equation}
  BF_{\bb{u}}^{-1}B^T \approx BM_{\bb{u}}^{-1}B^TF_p^{-1}M_p,
\end{equation}
which defines our approximate space-time pressure Schur complement \cref{eqn::XXfirstDef}.
\begin{remark}
The most recent edition of \cite{andyFIT} proposes to consider the commutator
between the discrete \emph{divergence} (rather than \emph{gradient}) operator,
and the reaction-advection-diffusion operator. The pressure Schur complement
approximation stemming from this choice is very similar to \cref{eqn::XXdef},
but with $M_p$ and $A_p$ swapped. As such, its application cost is virtually the same,
and preliminary numerical experiments did not provide significant improvements
over \cref{eqn::XXdef} in terms of convergence. Moreover, the form \cref{eqn::XXdef}
is more prone to simplifications similar to the one illustrated in \cref{eqn::XXsimp}
when solving Oseen equations, which is why we consider the gradient commutator \cref{eqn::XXdef}
for the remainder of this paper.
\end{remark}

\subsubsection{Comparison with single time-step case}
\label{sec::precon::STvsTS}
As briefly remarked at the end of \cref{sec::oseen}, the PCD preconditioner
\cite{PCDoriginal} which is used for the acceleration of the solution for the single
time-step system \cref{eqn::oseenBlock}, and which largely motivates our work,
for each instant $k$ is given by
\begin{equation}
  \mc{P}_{T,k}     \coloneqq\left[\begin{array}{cc}\mc{F}_{\bb{u},k}     &\mc{B}^T                                   \\&-\mc{X}_k     \end{array}\right]\quad\Longleftrightarrow\quad
  \mc{P}_{T,k}^{-1}\coloneqq\left[\begin{array}{cc}\mc{F}_{\bb{u},k}^{-1}&\mc{F}_{\bb{u},k}^{-1}\mc{B}^T\mc{X}_k^{-1}\\&-\mc{X}_k^{-1}\end{array}\right],
  \label{eqn::PTSingeStep}
\end{equation}
where the single time-step pressure Schur complement, $\mc{S}_{p,k}\coloneqq\mc{B}\mc{F}_{\bb{u},k}^{-1}\mc{B}^T$, is approximated as
\begin{equation}
  \mc{S}_{p,k}\approx\mc{X}_k \coloneqq \tilde{\mc{A}}_p \mc{F}_{p,k}^{-1}\mc{M}_p \quad\Longleftrightarrow\quad \mc{X}_k^{-1} \coloneqq \mc{M}_p^{-1} \mc{F}_{p,k} \tilde{\mc{A}}_p^{-1}
  , \qquad \text{with}\qquad k=1,\dots,N_t,
  \label{eqn::XXdefSingeStep}
\end{equation}
rendering it strikingly similar to \cref{eqn::XXdef}. The total computational
cost associated with applying the pressure Schur complement approximation
\cref{eqn::XXdefSingeStep} at each time-step is analogous to that of its
space-time counterpart; moreover, the application of \cref{eqn::XXdef} can be
conducted in parallel over time naturally. We can convince ourselves of this by considering the
block-diagonal structure of the matrices in \cref{eqn::XXdef}: applying the
inverses of $A_p$ and $M_p$ requires solving $N_t$ independent systems involving
$\tilde{\mc{A}}_p$ and $\mc{M}_p$, similarly to applying
\cref{eqn::XXdefSingeStep} for each of the $N_t$ time-steps. The sole overhead
with respect to sequential time-stepping comes from the \emph{bi-}diagonal
structure of $F_p$: the off-diagonal coupling between one time-step and the next
must be accounted for for \emph{each} application of the space-time
preconditioner, rather than once per time-step. This extra cost is however
negligible compared to the other operations.
In terms of set-up cost, the single time-step preconditioner is cheaper when solving
Stokes. In this case, in fact, the relevant operators in \cref{eqn::PTSingeStep}
are constant for each time-step, and they can hence be assembled once and for all and reused
throughout the time-stepping procedure. Conversely, in the space-time setting,
to ensure that the application of $P_T$ can be done in parallel, each processor
should assemble its own operators (at least in principle), thus
effectively multiplying the total set-up cost by $N_t$. Even the assembly operations,
however, can be conducted in parallel among the processors.
Moreover, this advantage is largely lost when considering time-dependent Oseen
or Navier-Stokes, as different convection-diffusion operators and the
corresponding preconditioners must be assembled for each time-step or
nonlinear iteration.

In light of these considerations, we conclude that the difference in
performance between using \cref{eqn::PTSingeStep} at each time-step, and
tackling the whole space-time system at once with \cref{eqn::PT}, is given by
two factors: (i) the difference in the total number of GMRES
iterations to convergence using the space-time block preconditioner
\cref{eqn::PT}, which we denote with $N_{it}^{ST}$, and the average number of
GMRES iterations to convergence per time-step using \cref{eqn::PTSingeStep},
identified by $N_{it}^0$; and (ii) the computational time associated with
inverting the velocity spatial operator $\mc{F}_{\bb{u},k}$ at each time-step
{versus} inverting the whole space-time velocity block $F_{\bb{u}}$. We
denote the latter as $C_{F_{\bb{u}}}$, while for the former we consider an
average of $C_{\mc{F}_{\bb{u}}}$ per time-step, so that the total time is
given by $C_{\mc{F}_{\bb{u}}} \cdot N_t$.
We expect our approach to be competitive if
\begin{equation}
  \frac{N_{it}^{ST}}{N_{it}^0}\,C_{F_{\bb{u}}} < C_{\mc{F}_{\bb{u}}}\,N_t.
  \label{eqn::STvsTS}
\end{equation}
The ratio $N_{it}^{ST} / N_{it}^0$ defines the overhead cost of the all-at-once
approach in terms of the total number of preconditioner applications. Numerical
results in \cref{sec::results::STvsSTS} indicate that this ratio is between $1$ and $2$
for most cases, in some extreme settings increasing to $3$. Thus, the overhead in
performing all-at-once block preconditioning in the space-time setting is quite
marginal compared to standard spatial preconditioning in sequential
time-stepping. The efficiency of our approach then directly depends on the speed-up that can be
obtained when solving a vector-based time-dependent (advection-)diffusion
problem using PinT or all-at-once techniques. Moreover, with a sufficiently
large number of MPI processes available, PinT and all-at-once solutions to
space-time (advection-)diffusion problems have shown significant speed-ups over
sequential time-stepping. For instance, one example in \cite{Falgout.2017}
demonstrated as large as a $325\times$ speed-up using the all-at-once multigrid
solution of space-time diffusion over sequential time-stepping. Thus,
coupling with the space-time block preconditioning approach developed here,
we expect a time-dependent Stokes problem with a comparable problem size in
space and time and the same number of MPI processes to lead to
a speed-up between $\sim100$-$300\times$ over sequential time-stepping.

\subsubsection{Defining the pressure operators}
\label{sec::precon::Apdef}
One issue associated with the definition of the pressure reaction--advection-diffusion
operator \cref{eqn::PRAD} and the pressure Laplacian \cref{eqn::Plaplacian}
involves the choice of BC to impose on said operators. In
\cite{PCDboundaryCond} and \cite[Chap.~9.2.2]{andyFIT}, it is suggested that
$\mc{F}_{p,k}$ should be assembled as if Robin BC,
\begin{equation}
  \mu\nabla p(\bb{s},t)\,\bb{n}(\bb{s}) - (\bb{w}(\bb{s},t)\cdot\bb{n}(\bb{s}))\,p(\bb{s},t) = 0,
\end{equation}
were imposed everywhere, while $\tilde{\mc{A}}_p$ should have Dirichlet BC on
the outflow boundary (that is, on $\bb{s}\in\partial\Omega :
\bb{u}(\bb{s})\cdot\bb{n}(\bb{s})>0$), and Neumann BC everywhere else. In
particular, this choice implies that $\mc{A}_p$ in \cref{eqn::PRAD} might differ
from $\tilde{\mc{A}}_p$. However, prescribing the same BC for the assembly of
both the pressure reaction-advection-diffusion operator and the pressure Laplacian has the
advantage of further simplifying \cref{eqn::XXdef} for the Stokes equations
\cref{eqn::stokes}. Under the assumption $\mc{A}_p\equiv\tilde{\mc{A}}_p$, in
fact, this becomes
\begin{equation}
X^{-1}
  \opif{=}{\left(\begin{aligned}\text{if }\mc{A}_p&\equiv\tilde{\mc{A}}_p\\\bb{w}&\equiv\bb{0}\end{aligned}\right)} \left[\begin{array}{ccc}
  \frac{\mc{A}_p^{-1}}{\Delta t}+\mu\mc{M}_p^{-1} &      & \\
 -\frac{\mc{A}_p^{-1}}{\Delta t}                  &\ddots& \\
                                                  &\ddots& \\
  \end{array}
  \right].
  \label{eqn::XXsimp}
\end{equation}
In comparison to \cref{eqn::XXdef}, we need to multiply by one fewer operator,
and inverting $\mc{M}_p$ and $\mc{A}_p$ can be done independently. Notice that a
preconditioner in this form is reminiscent of the Cahouet-Chabard preconditioner
\cite{CahouetChabard}. Keeping performance in mind, then, in our code we rather
choose to assemble both the pressure Laplacian and the pressure
reaction-advection-diffusion operator considering homogeneous Dirichlet BC at the outflow
of the domain, and homogeneous Neumann BC everywhere else, thus mimicking the
implementation in the IFISS package \cite{IFISS} (available at \cite{IFISSweb}).
Moreover, early experiments showed that the two approaches provide comparable
results in terms of acceleration of convergence for the test-cases considered.

We also remark that, using the BC prescribed above,
\cref{eqn::massStiffP::stiff} becomes singular when dealing with {fully
enclosed} Stokes flows, for which $\bb{u}(\bb{s})\cdot\bb{n}(\bb{s})=0$,
$\forall\bb{s}\in\partial\Omega$. In this context, in fact, the operator mimics
the action of a Laplacian with Neumann BC everywhere on the boundary, whose
kernel contains all constant functions. The GMRES procedure remains robust also
in this case \cite[Chap.~9.3.5]{andyFIT}, but we need to ensure that the chosen
solver for $\tilde{\mc{A}}_p$ can deal with singular matrices.

Finally, since with \cref{eqn::funcSpaceP} we only require $L^2$-regularity in
space for the pressure variable, the operator \cref{eqn::massStiffP::stiff}
might be ill-defined, depending on the discretisation chosen. A solution to
this issue when considering a discontinuous pressure field is detailed in
\cite[Chap.~9.2.1]{andyFIT}.

\subsection{Green's function approach}
\label{sec::precon::greenFun}
Alternatively to the small-commutator approach described in
\cref{sec::precon::smallComm}, another heuristic for recovering a good candidate
preconditioner for the Oseen equations has been proposed in \cite{tensorPC}. It
is based on mapping the operators appearing in \cref{eqn::oseen} to Fourier
space, finding the associated Green's tensor in Fourier domain, and taking
inspiration from its form to recover an approximation of the Schur complement.
Here we demonstrate that such principles can also be extended to the space-time
setting, and, in fact, arrive at the same Schur-complement approximation as
derived in \cref{eqn::XXdef}, providing additional validation for our choice of
preconditioner.

We start by arranging the operators appearing in \cref{eqn::oseen} in tensor
form,
\begin{equation}
  \left[\begin{array}{c|c}
          \mc{D}_\varphi & \nabla \TopSp\BotSp\\\hline
          \nabla\cdot    & 0\TopSp\BotSp
        \end{array}\right]
  \left[\begin{array}{c}
    \bb{u}(\bb{x},t)\TopSp\BotSp\\\hline
    p(\bb{x},t)     \TopSp\BotSp
  \end{array}\right]=\left[\begin{array}{c}
    \bb{f}(\bb{x},t)\TopSp\BotSp\\\hline
    0               \TopSp\BotSp
  \end{array}\right],
\end{equation}
where $\mc{D}_{(*)}$ denotes the block diagonal operator responsible for
applying $(*)$ to each component of a vector in $\mathbb{R}^d$, while $\varphi$
is defined as
\begin{equation}
  \varphi \coloneqq \frac{\partial }{\partial t} + \bb{w} \cdot \nabla - \mu \nabla^2.
  \label{eqn::OseenComponent}
\end{equation}
The (tensor) operator associated to \cref{eqn::oseen} is then given by
\begin{equation}
  \mc{L}\coloneqq \left[\begin{array}{c|c}
    \mc{D}_{\varphi} & \nabla \TopSp\BotSp\\\hline
    \nabla\cdot      & 0      \TopSp\BotSp
  \end{array}\right].
  \label{eqn::tensorOseen}
\end{equation}
We seek to recover a candidate approximation for the pressure Schur complement appearing in our preconditioner \cref{eqn::PT}, by following the form of the
Green's tensor $\mc{G}$ associated with \cref{eqn::tensorOseen}, under the
assumption of an advection field constant over the spatial domain
(but not necessarily in time), $\bb{w}(\bb{x},t)\equiv \bb{w}(t)$. To this purpose, we can
exploit the fact that the Fourier transform of the Green's tensor $\hat{\mc{G}}$
is given by the inverse of the Fourier transform of \cref{eqn::tensorOseen},
that is, $\hat{\mc{G}} = \hat{\mc{L}}^{-1}$.  The operator equivalent to
\cref{eqn::OseenComponent} in Fourier space is given by
\begin{equation}
  \hat{\varphi} = \frac{\partial }{\partial t} + \mathrm{i} \bb{w} \cdot \bb{k} + \mu |\bb{k}|^2,
  \label{eqn::OseenComponentFT}
\end{equation}
where $\bb{k}=\left[k_0,\dots,k_{d-1}\right]^T$ is the frequency vector, and $\mathrm{i}$
represents the imaginary unit; consequently, the Fourier equivalent of
\cref{eqn::tensorOseen} is given by
\begin{equation}
  \hat{\mc{L}} = \left[\begin{array}{c|c}
    \mc{D}_{\hat{\varphi}} & \mathrm{i}\bb{k} \TopSp\BotSp\\\hline
    \mathrm{i}\bb{k}\cdot           & 0       \TopSp\BotSp
  \end{array}\right].
  \label{eqn::tensorOseenFT}
\end{equation}
Let $\hat{\sigma}_p$ denote the (2,2) Schur complement in Fourier space,
\begin{equation}
  \hat{\sigma}_p \coloneqq - \mathrm{i}\bb{k}^T\mc{D}_{\hat{\varphi}^{-1}}\mathrm{i}\bb{k}
    = |\bb{k}|^2\hat{\varphi}^{-1}.
  \label{eqn::SchurCompFT}
\end{equation}
Then, using the classic formula for inversion of a $2\times2$ block operator based
on a block LDU decomposition, we recover
\begin{equation}
\begin{split}
\hat{\mc{G}}=\hat{\mc{L}}^{-1}
  &=\left[\begin{array}{c|c}
    \mc{D}_{\hat{\varphi}^{-1}}-\mc{D}_{\hat{\varphi}^{-1}}\bb{k}\hat{\sigma}_p^{-1}\bb{k}^T\mc{D}_{\hat{\varphi}^{-1}} & -\mc{D}_{\hat{\varphi}^{-1}}\mathrm{i}\bb{k}\hat{\sigma}_p^{-1}\TopSp\BotSp\\\hline
    -\hat{\sigma}_p^{-1} \mathrm{i}\bb{k}^T  \mc{D}_{\hat{\varphi}^{-1}}                                                         & \hat{\sigma}_p^{-1}\TopSp\BotSp
    \end{array}\right]\\
  &=\left(\hat{\varphi}|\bb{k}|^2\right)^{-1}\left[\begin{array}{c|c}
      \mc{D}_{|\bb{k}|^2}-\bb{k} \bb{k}^T & -\mathrm{i}\bb{k} \hat{\varphi} \TopSp\BotSp\\\hline
      -\mathrm{i}\bb{k}^T\hat{\varphi}             & \hat{\varphi}^2        \TopSp\BotSp
    \end{array}\right].
\label{eqn::GreenF}
\end{split}
\end{equation}
We focus our attention on the bottom-right block of \cref{eqn::GreenF}.
Transforming back from Fourier space to physical space, and using the properties
of Fourier transforms, namely the correspondence of the operators
$|\bb{k}|^2\leftrightarrow-\nabla^2$ and
$-\bb{k}\bb{k}^T\leftrightarrow\nabla(\nabla\cdot)$, as well as the identity
$\nabla\times(\nabla\times) = \nabla(\nabla\cdot) -  \nabla^2$, we get an
expression for the Green's tensor:
\begin{equation}
\mc{G} = \left[\begin{array}{c|c}
    (\nabla\times(\nabla\times))\mc{D}_{\mc{G}_{\varphi\nabla^2}} & (-\nabla)\mc{G}_{-\nabla^2}\TopSp\BotSp\\\hline
    (-\nabla\cdot)\mc{D}_{\mc{G}_{-\nabla^2}}                     & (\frac{\partial}{\partial t}+\bb{w}\cdot\nabla-\mu\nabla^2)\mc{G}_{-\nabla^2}\TopSp\BotSp
  \end{array}\right],
\label{eqn::Green}
\end{equation}
where $\mc{G}_{-\nabla^2}$ is the fundamental solution for the negative
Laplacian operator, while $\mc{G}_{\varphi\nabla^2}$ is the fundamental solution
for the operator $\nabla^2\left(\frac{\partial }{\partial t} + \bb{w} \cdot
\nabla - \mu \nabla^2\right)$. Keeping only the bottom-right block and
discarding the rest, we get that the space-time pressure Schur complement
operator $\sigma_p$ satisfies
\begin{equation}
  \sigma_p^{-1}q \approx \left(\left(\frac{\partial }{\partial t} + \bb{w} \cdot \nabla - \mu \nabla^2\right)\mc{G}_{-\nabla^2}\right)*q, \qquad\qquad\forall q\in Q(\Omega),
\end{equation}
and on this we base its approximation. Considering that
$\mc{G}_{-\nabla^2}*q$ should represent the solution to a (negative) Laplacian
defined on the pressure space, we can discretise it using $\tilde{\mc{A}}_p^{-1}
q$ defined in \cref{eqn::Plaplacian}. Analogously, the space-time operator
acting on the pressure space can be discretised via $F_p \approx \frac{\partial
}{\partial t} + \bb{w}\cdot \nabla - \mu \nabla^2$, introduced in
\cref{eqn::STFp}.

Altogether, this analysis suggests that a discrete approximation to the
space-time pressure Schur complement $S_p$ should satisfy
\begin{equation}
  S_p^{-1} \approx X^{-1} \coloneqq M_p^{-1}F_p A_p^{-1},
  \label{eqn::XXdefFourier}
\end{equation}
where $M_p$ and $A_p$ are defined as in \cref{eqn::MpdiagApdiag}. Note, we
re-scale the operator by the pressure mass matrix for similar reasons as
discussed for \cref{eqn::commutator}, and as has also been applied to the single
time-step case in \cite{tensorPC}. We can see that the Green's function
heuristic followed in this section yields the same approximation to the
space-time pressure Schur complement as found via the small-commutator approach
in \cref{sec::precon::smallComm}. This is not necessarily surprising, as both
analyses in some sense rely on neglecting boundary conditions, but the fact that
two approaches provide the same approximation \cref{eqn::XXdefFourier} gives
confidence in the validity of the proposed preconditioner.


\subsection{Eigenvalues clustering}
\label{sec::precon::eigs}
To provide an indication of the effectiveness of \cref{eqn::PT} as a
preconditioner, we investigate how the eigenvalues of the resulting
preconditioned system spread in the complex plane. Although eigenvalues
are not necessarily indicative of fast convergence for GMRES, particularly
for highly nonsymmetric systems \cite{convGMRES}, poor clustering of
eigenvalues is likely to yield poor convergence, while nicely bounded
and clustered eigenvalues is still often indicative of an effective 
preconditioning procedure.

We are interested in solving the following generalised eigenvalue problem:
\begin{equation}
  \left[\begin{array}{cc}
    F_{\bb{u}} & B^T \\
    B          &
  \end{array}\right]\left[\begin{array}{c}
    \bb{u} \\
    \bb{p} 
  \end{array}\right]=\lambda\left[\begin{array}{cc}
    F_{\bb{u}} & B^T\\
               & -X
  \end{array}\right]\left[\begin{array}{c}
    \bb{u} \\
    \bb{p}
  \end{array}\right].
  \label{eqn::GEPbTprec}
\end{equation}
The system presents $\lambda=1$ as a solution with multiplicity $N_\bb{u}$; the
remaining ones are given by satisfying
\begin{equation}
  \left(X^{-1}BF_{\bb{u}}^{-1}B^T - \lambda I\right)\bb{p} = 0.\\
\end{equation}
From this, it can be seen that the general eigenvalues of \cref{eqn::GEPbTprec}
directly relate to the eigenvalues of the matrix $X^{-1}BF_{\bb{u}}^{-1}B^T$.
Substituting our approximation to the space-time pressure Schur complement
\cref{eqn::XXdef}, we have that the operator of interest becomes
\begin{equation}
  X^{-1}BF_{\bb{u}}^{-1}B^T = M_p^{-1}F_p A_p^{-1}BF_{\bb{u}}^{-1}B^T.
  \label{eqn::GEprob}
\end{equation}
This matrix is block lower triangular because $F_p$ and $F_{\bb{u}}^{-1}$ are
both block lower triangular and the remaining factors are block diagonal. Thus, to find its
eigenvalues, we just need to recover the eigenvalues of the blocks on the main
diagonal. These blocks consist of the operators
\begin{equation}
  \mc{M}_p^{-1}\mc{F}_{p,k}\mc{A}_p^{-1}\mc{B}\mc{F}_{\bb{u},k}^{-1}\mc{B}^T, \qquad\text{with}\qquad k=1,\dots,N_t.
  \label{eqn::blockGEP}
\end{equation}
Although we were unable to derive a tight theoretical bound on the eigenvalues
of \cref{eqn::blockGEP}, we can demonstrate nice clustering of preconditioned
eigenvalues numerically. In \cref{fig::eigsClust} are shown plots of eigenvalues
of \cref{eqn::blockGEP} in the complex plane, for varying values of $\Delta t$, $\Delta x$, and intensity of
the advection field $\bb{w}(\bb{x},t)$ in \cref{eqn::oseen} (the latter is
regulated by the P\'eclet number, Pe; see \cref{pb::doubleGlazing} for its
definition). We can see that in all cases the eigenvalues are nicely bounded,
they remain away from the origin, and are largely independent of variations in
the discretisation parameters, providing semi-rigorous theoretical support for
the proposed preconditioner.

\pgfplotsset{
  eigsPlotStyle/.style={
    tick label style={font=\small},
    xlabel= $\Re(\lambda)$,
    ylabel= $\Im(\lambda)$,
    axis background/.style={fill=gray!10},
    cycle list = { {{YlGnBu-L},only marks,mark=+}, {{YlGnBu-G},only marks,mark=o}, {{YlGnBu-E},only marks,mark=triangle*}},
    xmin = -0.2, xmax = 2.8, ymin = -0.7, ymax = 0.7,
    xmajorgrids= true,
    legend = north east,
    legend cell align={left},
    legend style={font=\small},
    filter discard warning=false,
    x filter/.expression={x<1e-14 ? nan : x},
  }
}
\begin{figure}[!t]
\centering
\begin{minipage}{0.49\textwidth}
  \centering
  \begin{tikzpicture}[baseline, scale=0.7]
      \begin{axis}[eigsPlotStyle]   
        \addplot table[scatter,x index=2,y index=3]{data/Pb4_Prec1_STsolve0_oU2_oP1_Pe100.000000/eigs_r4.dat}; \addlegendentry{$\text{Pe}=100$};
        \addplot table[scatter,x index=2,y index=3]{data/Pb4_Prec1_STsolve0_oU2_oP1_Pe10.000000/eigs_r4.dat};  \addlegendentry{$\text{Pe}=10$};
        \addplot table[scatter,x index=2,y index=3]{data/Pb1_Prec1_STsolve0_oU2_oP1/eigs_r4.dat};              \addlegendentry{$\text{Pe}=0$};
        \node[anchor=south west, fill=white!10 ] at (axis cs: -0.2,-0.7) {$\begin{array}{l}\Delta x\approx2^{-4}\\\Delta t=2^{-3}\end{array}$};
      \end{axis}
  \end{tikzpicture}
\end{minipage}
\hfill
\begin{minipage}{0.49\textwidth}
  \centering
  \begin{tikzpicture}[baseline, scale=0.7]
      \begin{axis}[eigsPlotStyle]   
        \addplot table[scatter,x index=6,y index=7]{data/Pb4_Prec1_STsolve0_oU2_oP1_Pe100.000000/eigs_r4.dat}; \addlegendentry{$\text{Pe}=100$};
        \addplot table[scatter,x index=6,y index=7]{data/Pb4_Prec1_STsolve0_oU2_oP1_Pe10.000000/eigs_r4.dat};  \addlegendentry{$\text{Pe}=10$};
        \addplot table[scatter,x index=6,y index=7]{data/Pb1_Prec1_STsolve0_oU2_oP1/eigs_r4.dat};              \addlegendentry{$\text{Pe}=0$};
        \node[anchor=south west, fill=white!10 ] at (axis cs: -0.2,-0.7) {$\begin{array}{l}\Delta x\approx2^{-4}\\\Delta t=2^{-5}\end{array}$};
      \end{axis}
  \end{tikzpicture}
\end{minipage}
\par
\begin{minipage}{0.49\textwidth}
  \centering
  \begin{tikzpicture}[baseline, scale=0.7]
      \begin{axis}[eigsPlotStyle]   
        \addplot table[scatter,x index=2,y index=3]{data/Pb4_Prec1_STsolve0_oU2_oP1_Pe100.000000/eigs_r6.dat}; \addlegendentry{$\text{Pe}=100$};
        \addplot table[scatter,x index=2,y index=3]{data/Pb4_Prec1_STsolve0_oU2_oP1_Pe10.000000/eigs_r6.dat};  \addlegendentry{$\text{Pe}=10$};
        \addplot table[scatter,x index=2,y index=3]{data/Pb1_Prec1_STsolve0_oU2_oP1/eigs_r6.dat};              \addlegendentry{$\text{Pe}=0$};
        \node[anchor=south west, fill=white!10 ] at (axis cs: -0.2,-0.7) {$\begin{array}{l}\Delta x\approx2^{-6}\\\Delta t=2^{-3}\end{array}$};
      \end{axis}
  \end{tikzpicture}
\end{minipage}
\hfill
\begin{minipage}{0.49\textwidth}
  \centering
  \begin{tikzpicture}[baseline, scale=0.7]
      \begin{axis}[eigsPlotStyle]   
        \addplot table[scatter,x index=6,y index=7]{data/Pb4_Prec1_STsolve0_oU2_oP1_Pe100.000000/eigs_r6.dat}; \addlegendentry{$\text{Pe}=100$};
        \addplot table[scatter,x index=6,y index=7]{data/Pb4_Prec1_STsolve0_oU2_oP1_Pe10.000000/eigs_r6.dat};  \addlegendentry{$\text{Pe}=10$};
        \addplot table[scatter,x index=6,y index=7]{data/Pb1_Prec1_STsolve0_oU2_oP1/eigs_r6.dat};              \addlegendentry{$\text{Pe}=0$};
        \node[anchor=south west, fill=white!10 ] at (axis cs: -0.2,-0.7) {$\begin{array}{l}\Delta x\approx2^{-6}\\\Delta t=2^{-5}\end{array}$};
      \end{axis}
  \end{tikzpicture}
\end{minipage}
\par
\caption{\label{fig::eigsClust} Eigenvalue distribution of the operators
\cref{eqn::blockGEP}, corresponding to \cref{pb::doubleGlazing}, for varying
values of $\Delta t$, $\Delta x$, and P\'eclet number, $\text{Pe}$, which
represents the intensity of the advection field. Here, Pe denotes the \emph{continuous}
P\'eclet number (i.e., not scaled by mesh spacing); see \cref{sec::modelProblems} for its
definition (notice choosing $\text{Pe}=0$ corresponds to \cref{pb::drivenCavity}).
The eigenvalues are computed numerically using the MATLAB function \texttt{eigs}.}
\end{figure}
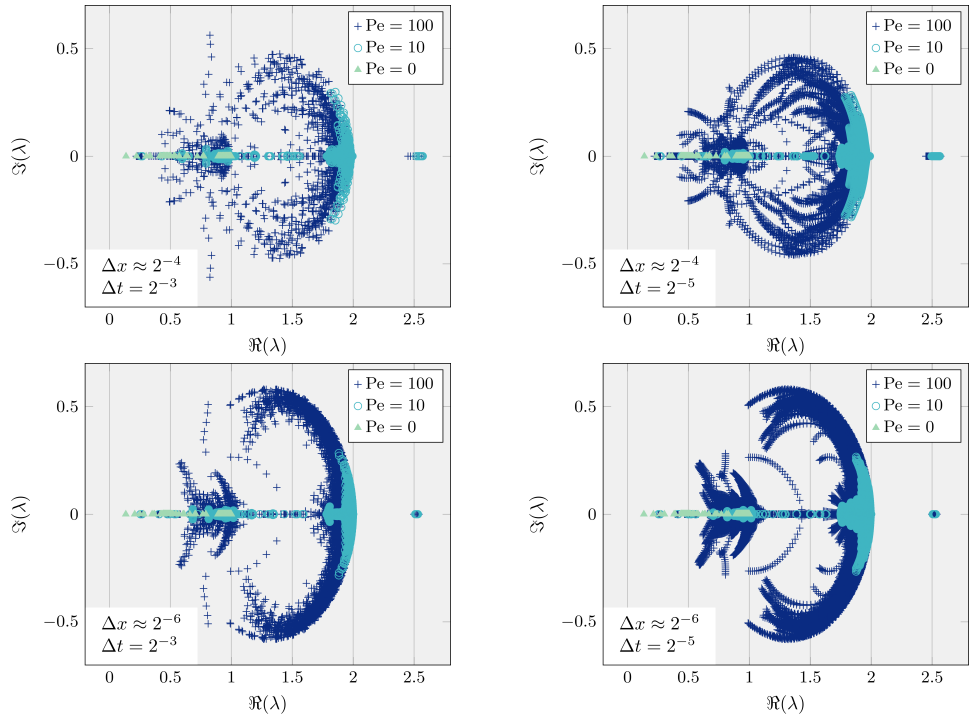

It is worth pointing out that increasing the P\'eclet number has
the effect of pushing some eigenvalues towards $2$, and spreading more of them
away from the real axis. It also makes the distributions more
sensitive with respect to the level of refinement in the spatial discretisation,
as can be seen comparing top and bottom plots in \cref{fig::eigsClust}: the dark
blue tokens become more orderly distributed as $\Delta x$ decreases. These
considerations makes us wary of a possible degradation in the performance of the
preconditioner as advection becomes predominant, an issue that is observed and
discussed in \cref{sec::results::advection}. In general, however,
most of the eigenvalues remain clustered inside $[0.1,2] \times [-0.6,0.6]$ in
the complex plane.

\end{section}
\begin{section}{Results}
\label{sec::results}
To demonstrate the performance of the preconditioner introduced in
\cref{sec::precon}, we consider its application for the solution of a variety of
flow configurations. For simplicity, we focus on problems defined on
$\Omega\subset\mathbb{R}^2$, although we point out that the applicability of
\cref{eqn::PT} is by no means limited to 2D simulations. The test-cases hereby
proposed have been adapted from problems used extensively in the literature (see
for example \cite[Chap.~3.1 and Chap.~6.1]{andyFIT}), and are briefly described
in \cref{sec::modelProblems}. The experiments conducted provide evidence for the
optimal scalability of \cref{eqn::PT}, as well as its potential for speed-up via
time-parallelisation. Additionally, we investigate simple extensions to more
general frameworks, such as the solution of nonlinear incompressible flow, and its
applicability in advection-dominated regimes.

\subsection{Model problems}
The first three problems described in this section are test-cases for the time-dependent Stokes equations \cref{eqn::stokes}, while the last one is
designed for Oseen \cref{eqn::oseen}. They are defined as follows.
\label{sec::modelProblems}
\begin{problems}
	\item\emph{Driven cavity flow.}\label{pb::drivenCavity}
  This represents a fully enclosed flow, defined on a square spatial domain
  $\bb{x}=(x,y)^T\in [0,1]^2 \eqqcolon \Omega_{\square}$. No forcing term is
  considered, and homogeneous Dirichlet BC are imposed on the bottom, left and
  right sides of the domain. To accelerate the flow, we prescribe the velocity
  on the top side
  \begin{equation}
  	\bb{u}_C(\bb{x}|_{y=1}, t) = 8tx(1-x)(2x^2-2x+1).
	\end{equation}
  This profile is regularised so as to smoothly match the BC on the left and right
  sides. Notice that, in order to introduce time-dependency, the velocity is
  ramped up linearly from a quiet state: this is done in the following
  test-cases as well. An example solution for this problem is shown in
  \cref{fig::solCavity}.
  \begin{figure}[!t]
    \begin{minipage}{0.45\textwidth}
      \centering
      \includegraphics[trim={2.42cm 0.605cm 2.42cm 0.605cm},clip,width=0.9\textwidth]{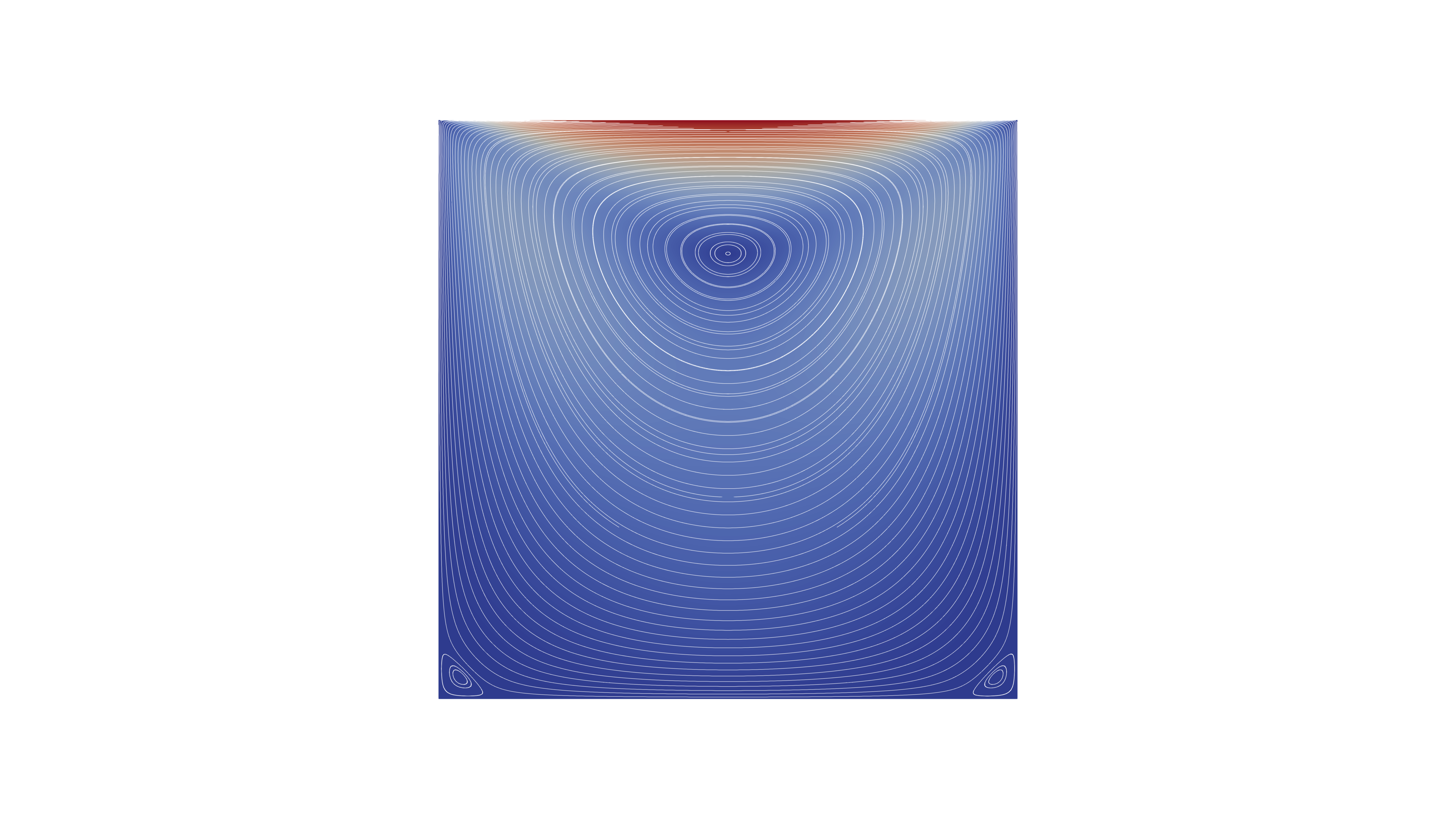}
    \end{minipage}
    \hfill
    \begin{minipage}{0.45\textwidth}
      \centering
      \includegraphics[trim={2.42cm 0.605cm 2.42cm 0.605cm},clip,width=0.9\textwidth]{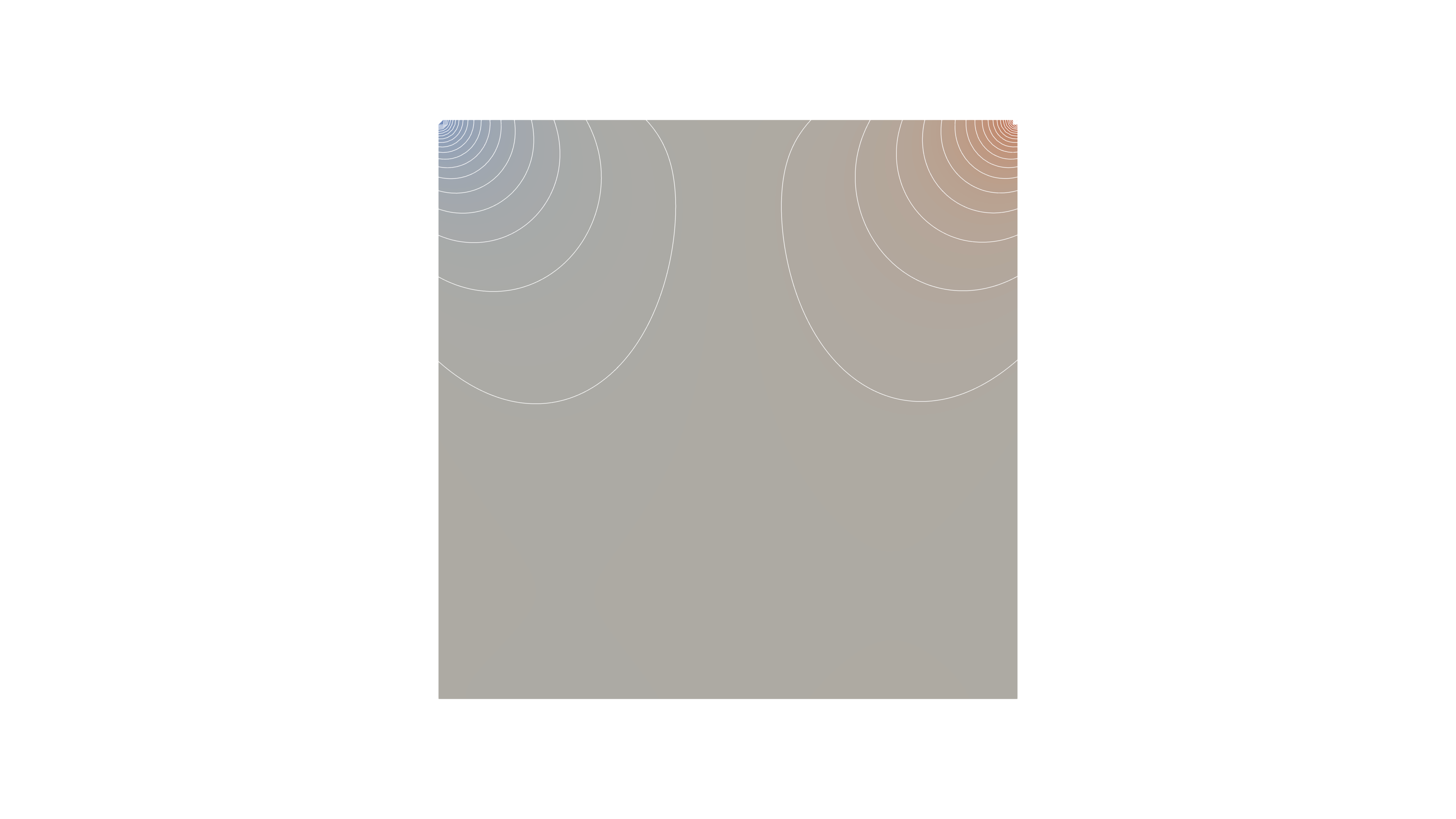}
    \end{minipage}
    \caption[Solution for driven cavity flow]{Velocity magnitude and streamlines (left), and pressure contour plot (right) for an example solution of the driven cavity flow \cref{pb::drivenCavity}.}
    \label{fig::solCavity}
  \end{figure}

	\item\emph{Poiseuille flow.}\label{pb::poiseuille}
  Also defined on $\Omega_{\square}$, this set-up models laminar flow along a
  channel, for which the following analytical solution can be recovered:
  \begin{equation}
	  \bb{u}_P(\bb{x},t) = 4t\left[\begin{array}{c}y(1-y)\\0\end{array}\right],\quad p_P(\bb{x},t) = 8(1-x).
	\end{equation}
  \begin{figure}[!b]
    \begin{minipage}{0.45\textwidth}
      \centering
      \includegraphics[trim={2.42cm 0.605cm 2.42cm 0.605cm},clip,width=0.9\textwidth]{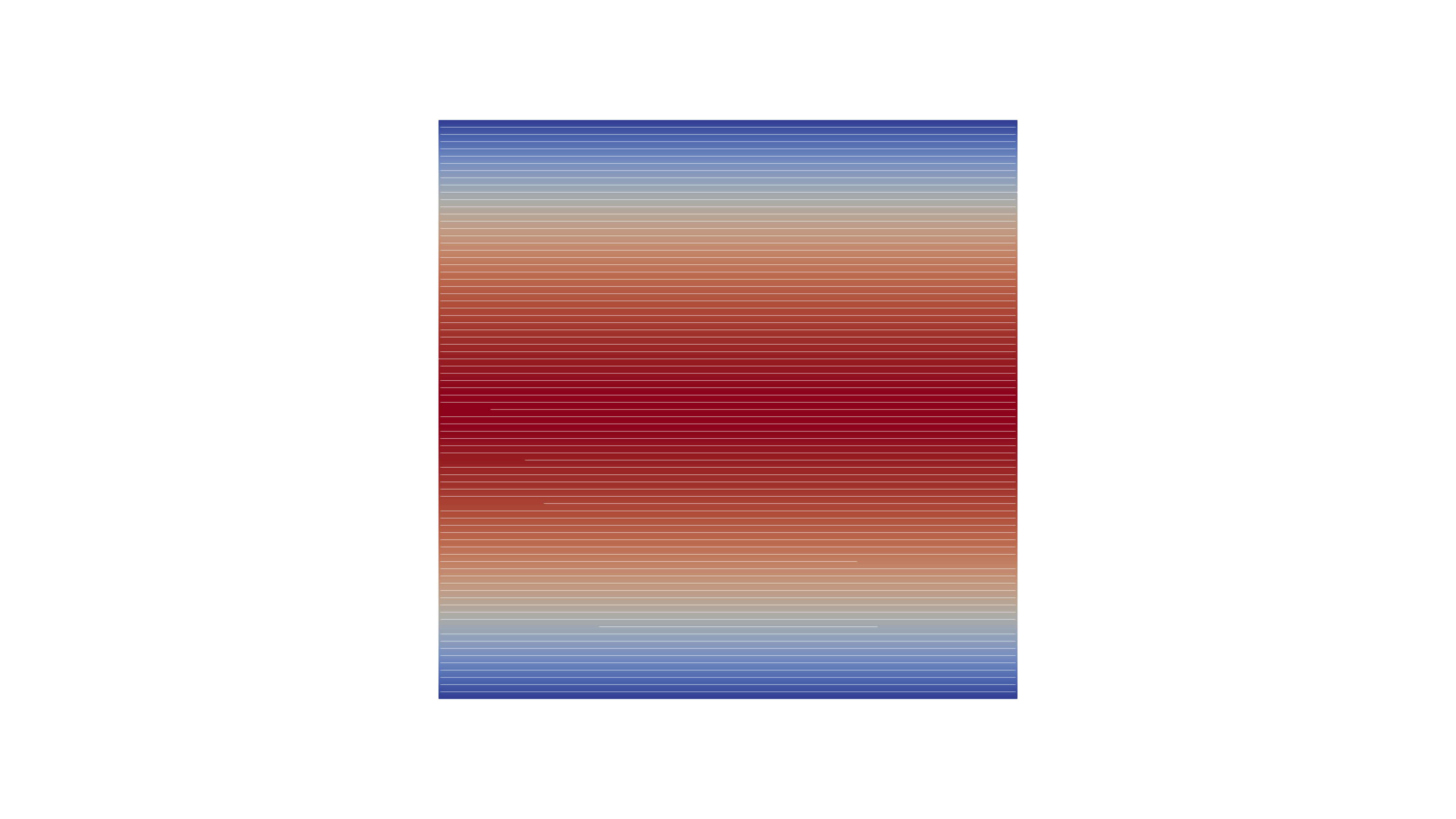}
    \end{minipage}
    \hfill
    \begin{minipage}{0.45\textwidth}
      \centering
      \includegraphics[trim={2.42cm 0.605cm 2.42cm 0.605cm},clip,width=0.9\textwidth]{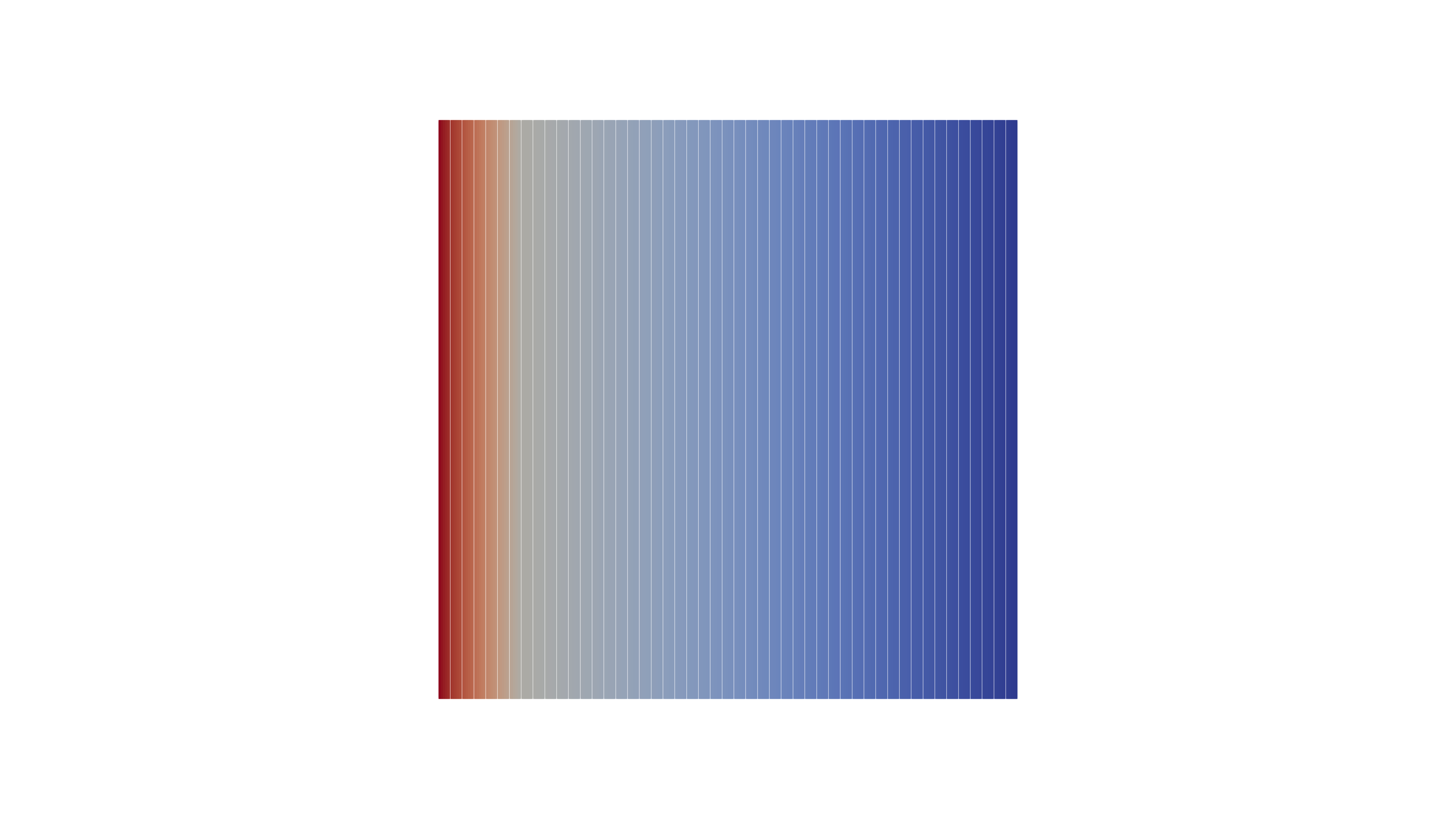}
    \end{minipage}
    \caption[Solution for Poiseuille flow]{Velocity magnitude and streamlines (left), and pressure contour plot (right) for the analytical solution of the Poiseuille flow \cref{pb::poiseuille}.}
    \label{fig::solPoiseuille}
  \end{figure}
  The forcing term is defined so as to counterbalance the time-derivative:
  $\boldsymbol{f}_P = 4[y(1-y),0]^T$. No-slip conditions are included in the top
  and bottom sides of the channel, and a parabolic inflow profile is prescribed
  at $x=0$. At the outflow $x=1$, homogeneous Neumann BC are imposed. The
  analytical solution is plotted in \cref{fig::solPoiseuille}.

	\item\emph{Flow over backward-facing step.}\label{pb::backStep}
  This is another example of flow through a channel, except the fluid undergoes
  a sudden expansion. The domain is a stretched L-shape,
  $\Omega_L=[0,8]\times[0,1]\cup[1,8]\times[0,-1]$, and the BC imposed are
  similar to those for Poiseuille (the inflow, where the parabolic velocity
  profile is prescribed, is on the side $x=0$). As for \cref{pb::drivenCavity},
  no forcing term is included; an example solution is provided in
  \cref{fig::solStep}.
  \begin{figure}[!t]
    \centering
    \includegraphics[trim={1.82cm 1.525cm 1.82cm 1.525cm},clip,width=0.9\textwidth]{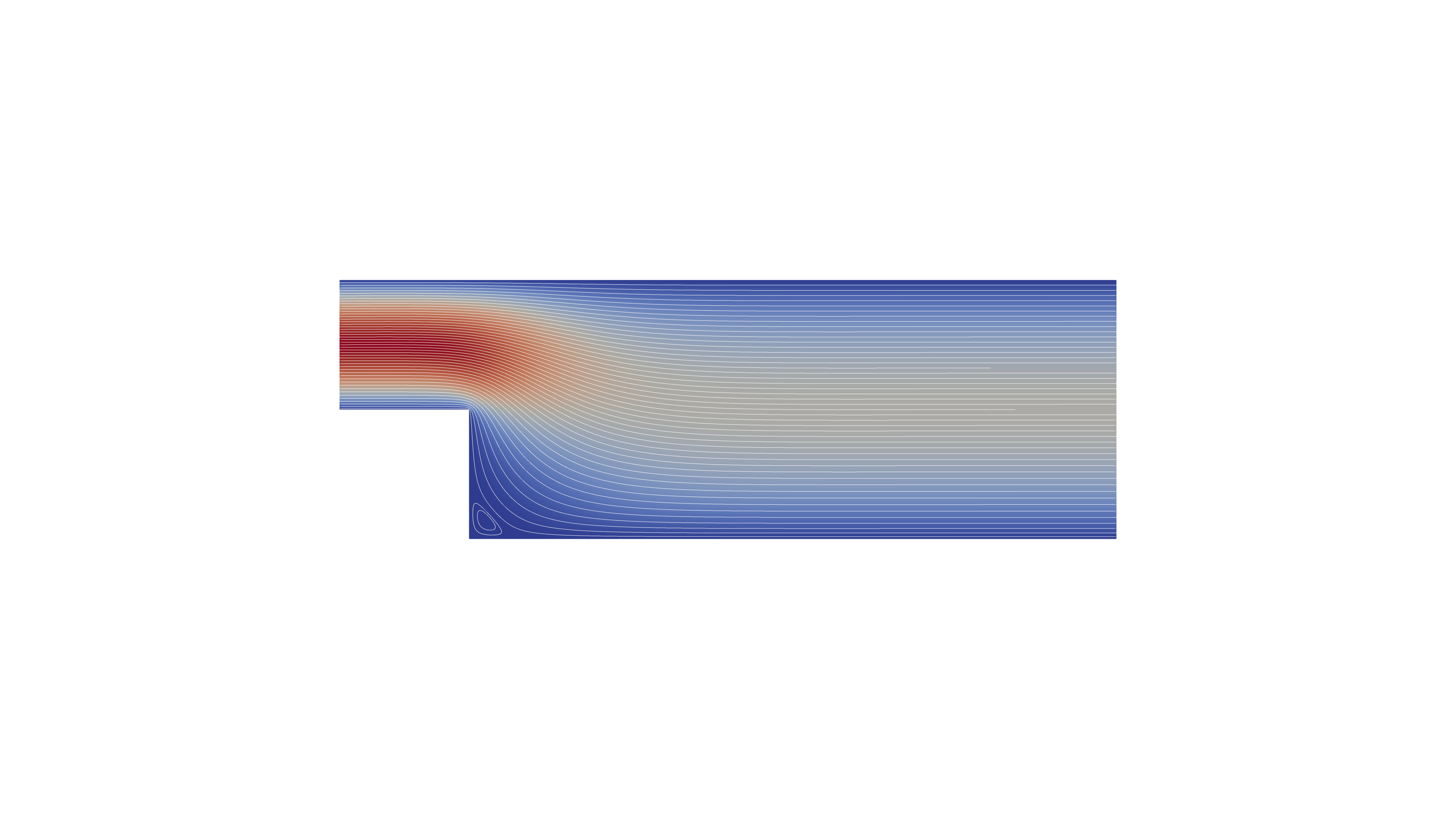}\\
    \includegraphics[trim={1.82cm 1.525cm 1.73cm 1.525cm},clip,width=0.9\textwidth]{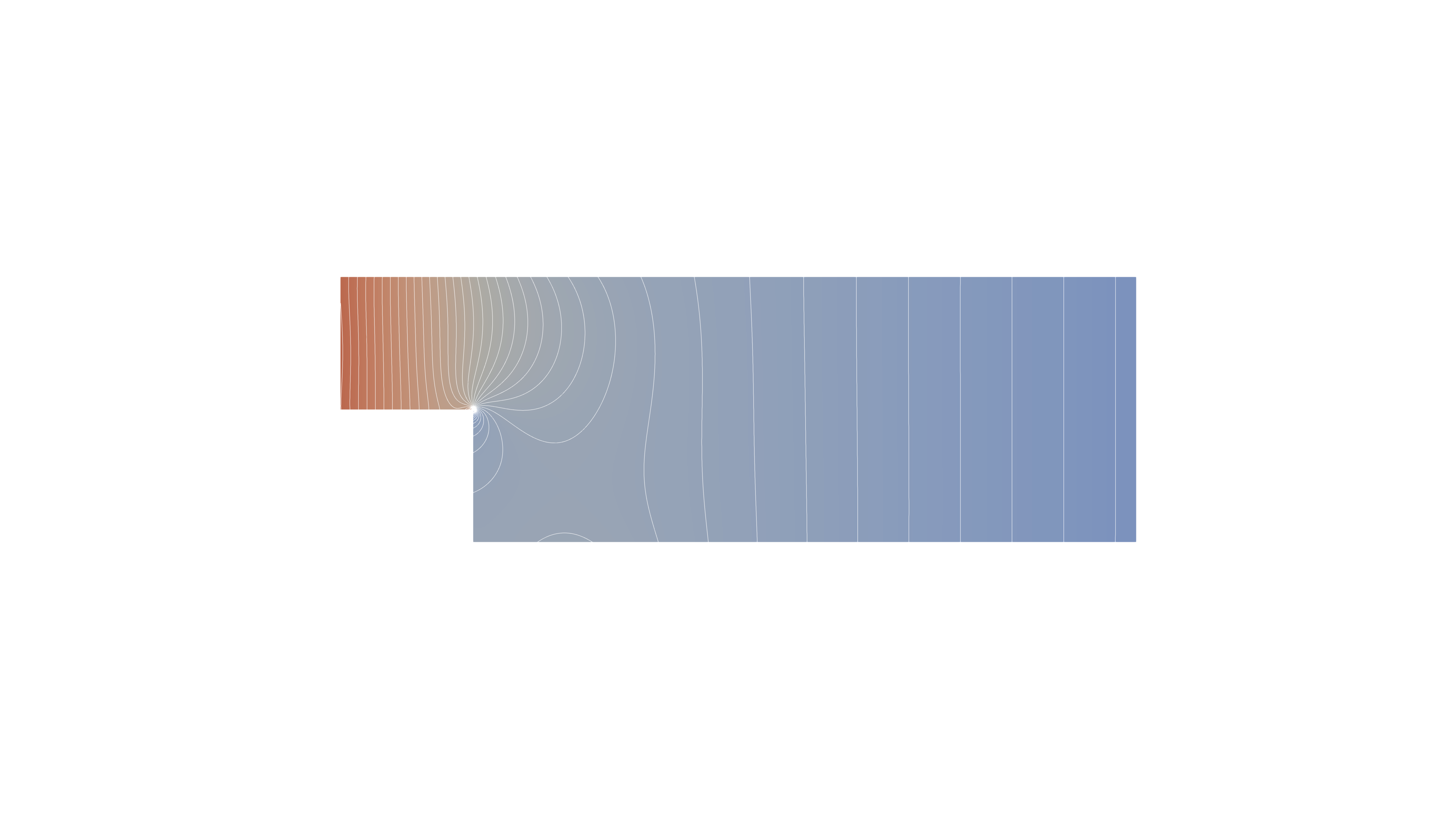}
    \caption[Solution for backward-facing step]{Velocity magnitude and streamlines (top), and pressure contour plot (bottom) for an example solution of the flow over backward-facing step \cref{pb::backStep}.}
    \label{fig::solStep}
  \end{figure}

	\item\emph{Double-glazing.}\label{pb::doubleGlazing}
	The set-up for this problem is analogous to that of the driven
  cavity flow, except now the effect of a recirculating wind is included in the
  system:
  \begin{equation}
  	\bb{w}_G(\bb{x},t) = 2t\left[\begin{array}{r}-(2y-1)(4x^2-4x+1)\\(2x-1)(4y^2-4y+1)\end{array}\right]\mu\text{Pe},
  \end{equation}
  whose intensity can be adjusted by tweaking the value of the P\'eclet number
  $\text{Pe}$ \cite[Chap.~13.2]{quarterNM}.
  This parameter describes the relative intensity of advective transport with
  respect to diffusive transport, and for a generic advection field $\bb{w}(\bb{x},t)$ it can
  be defined as
  \begin{equation}
    \text{Pe}\coloneqq \frac{LU}{\mu} = \frac{L}{\mu}\left(\frac{1}{T-T_0}\int_{T_0}^T \max_{\bb{x}\in\Omega}\bb{w}(\bb{x},t)\,dt\right),
    \label{eqn::peclet}
  \end{equation}
  where $L$ is a characteristic length of the spatial domain ($L=1$ for this problem), and $U$ a characteristic speed (considered averaged over the temporal domain).
  The impact that the added advection field has on the solution is illustrated in \cref{fig::solGlazing}.
  \begin{figure}[!t]
    \begin{minipage}{0.45\textwidth}
      \centering
      \includegraphics[trim={2.42cm 0.605cm 2.42cm 0.605cm},clip,width=0.9\textwidth]{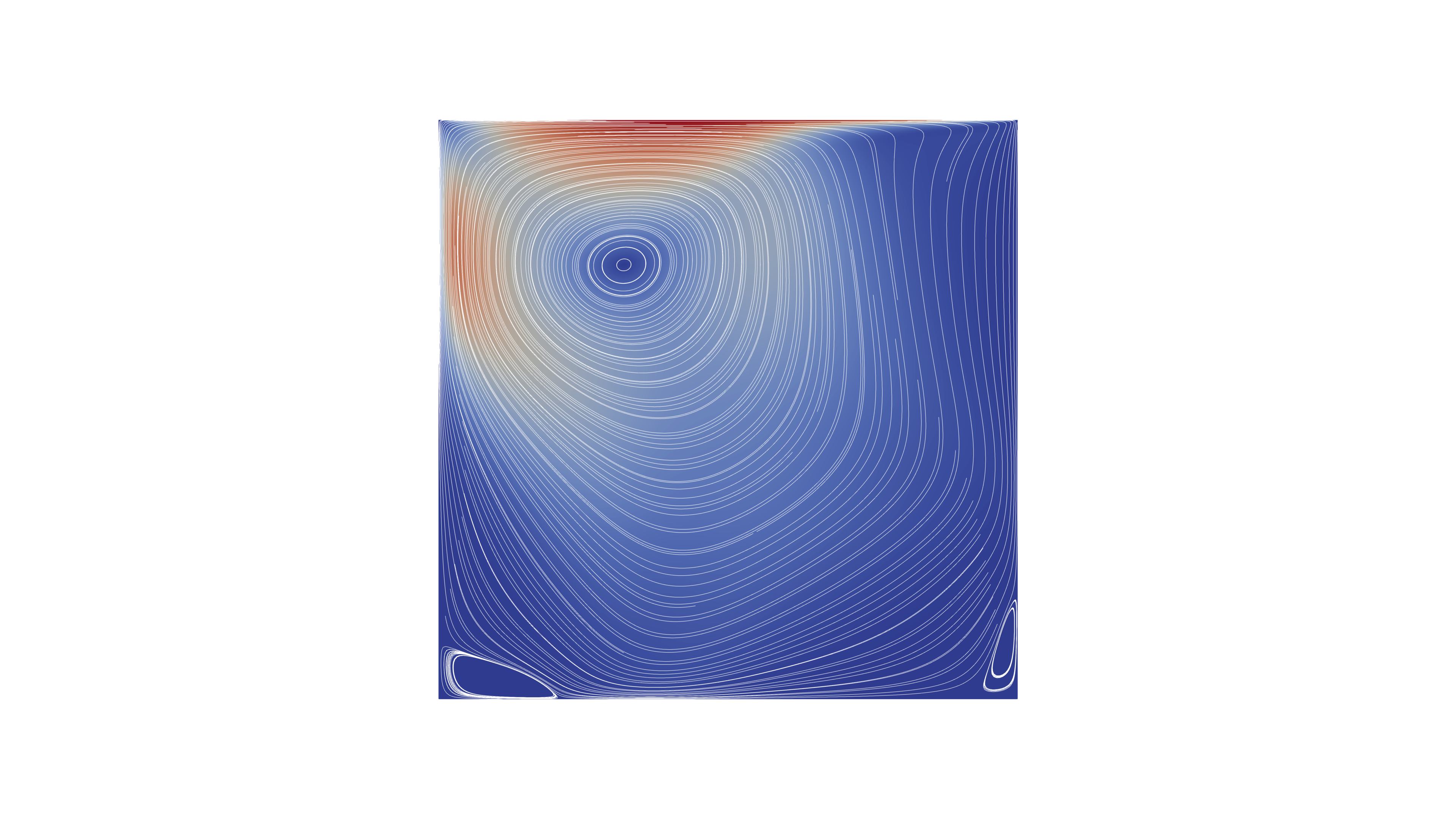}
    \end{minipage}
    \hfill
    \begin{minipage}{0.45\textwidth}
      \centering
      \includegraphics[trim={2.42cm 0.605cm 2.42cm 0.605cm},clip,width=0.9\textwidth]{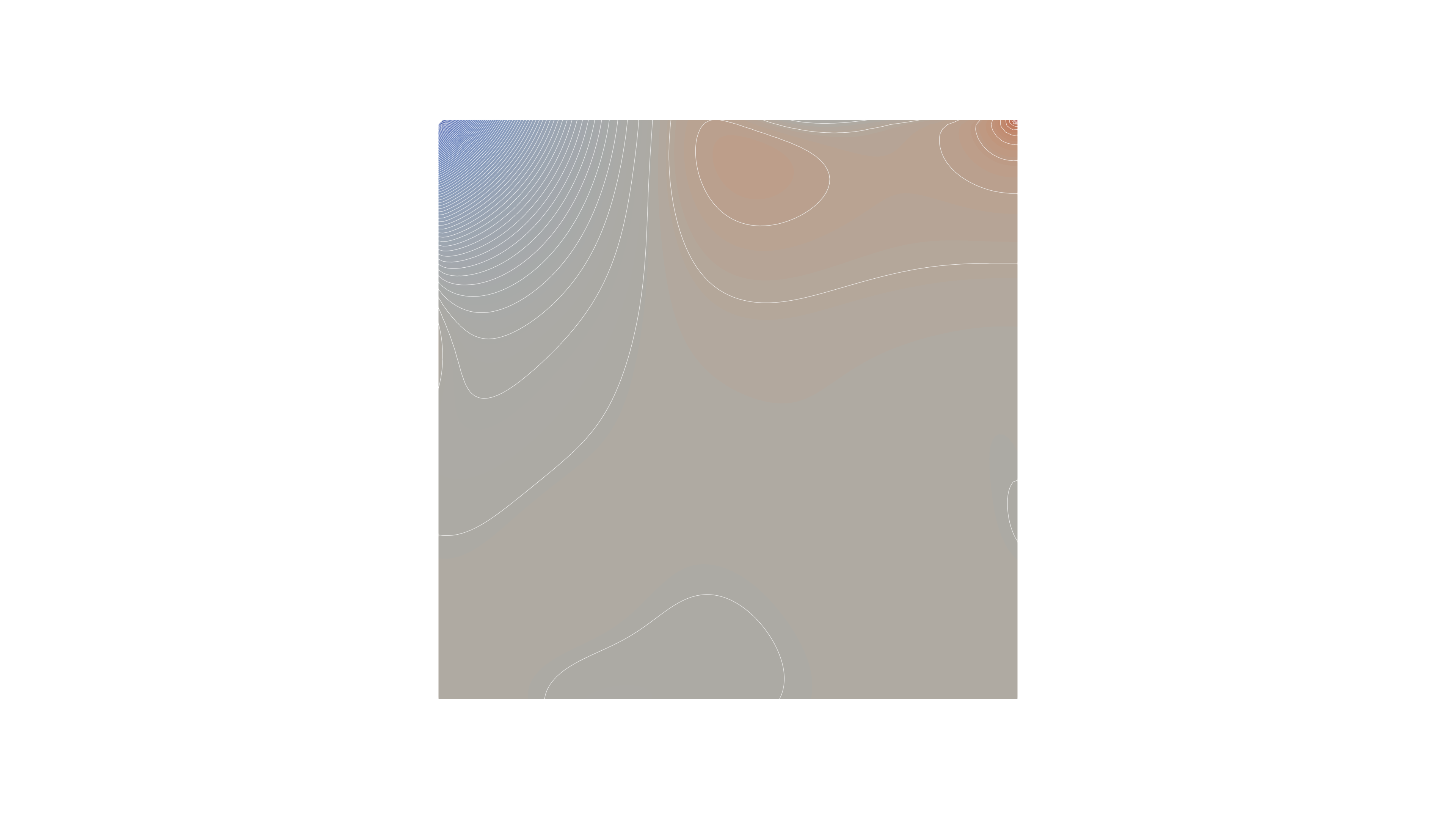}
    \end{minipage}
    \caption[Solution for double-glazing flow]{Velocity magnitude and streamlines (left), and pressure contour plot (right) for an example solution of the double-glazing flow \cref{pb::doubleGlazing} with $\text{Pe}=10$.}
    \label{fig::solGlazing}
  \end{figure}
\end{problems}

For each test-case the temporal domain is given by $t\in[0,1]$, and the
viscosity parameter is fixed to $\mu=1$. The discrete functional spaces are
approximated using Taylor-Hood finite elements: $P_2$ for velocity and $P_1$ for
pressure \cite[Chap.~17.4]{quarterNM}, both defined on a triangular mesh.
The code developed for our experiments is publicly available at the repository
\cite{ownRepo}. It makes use of the MFEM finite element package \cite{MFEM} for the
assembly of the finite element matrices, and we rely on the PETSc \cite{PETSC}
and \textit{hypre} \cite{HYPRE} implementations of the various solvers used. Unless
otherwise specified, the space-time initial guess for the GMRES iterations is null (apart
from the Dirichlet nodes, where we substitute the given data),
although preliminary experiments with a random initial guess
(uniformly distributed in $\left[-1,1\right]$) show negligible sensitivity of the convergence
behaviour to this choice.

\subsection{Performance of preconditioner}
In this section, we conduct a series of experiments aimed at measuring the
effectiveness of the proposed preconditioner \cref{eqn::PT} in accelerating GMRES.
We chose the total number of iterations to convergence as
the main measure of performance. The reason behind this preference (rather than,
for example, computational time), lies in the fact that we believe it provides a
more indicative metric for assessing the efficacy of our preconditioning
strategy. This is in fact less dependent on the details of the actual
implementation of the underlying solvers, and particularly of the one for the
space-time velocity block \cref{eqn::STFu}, the optimal design of which is
beyond the purpose of this manuscript: indeed a key point of our approach is the
flexibility it provides in choosing such a solver.

\subsubsection{Exact {versus} approximate solvers}
As a first experiment, we test directly the application of our preconditioner to
the solution of the model problems introduced in \cref{sec::modelProblems}. We
provide two different set-ups for the solvers involved in the application of
\cref{eqn::PT}: an \emph{ideal} one, in which we make use of {exact}
solvers, in order to provide a {best-case} scenario for the performance of
such preconditioner; and an \emph{approximate} case, which instead gives us a
measure of the performance degradation we can expect by applying iterative
solvers instead of exact solvers, which is more realistic in practice.

\begin{table}[t!]\footnotesize
\caption{\label{tab::PconvIdealVSapprox}
Number of GMRES iterations to $10^{-10}$ residual tolerance, when
applied to the monolithic space-time system \cref{eqn::STdisc}
right-preconditioned with $P_T^{-1}$ \cref{eqn::PT}. Different levels of refinement are
considered both for the spatial (rows) and temporal meshes (columns). The four
sets of results refer to the four test problems introduced in
\cref{sec::modelProblems} (\cref{pb::doubleGlazing} has P\'eclet number
$\text{Pe}=10$). Results on the left of each column are recovered using exact
solvers, and results on the right of each column (in brackets) are recovered using approximate
iterative solvers (inversion of $\mc{M}_p$ is approximated with 8 Jacobi-preconditioned
Chebyshev iterations, inversion of $\mc{A}_p$ is approximated using 15
iterations of AMG, and the solution of the space-time velocity block
$F_{\bb{u}}$ is recovered applying 15 iterations of GMRES preconditioned by
AIR). The outer solver in the latter case is FGMRES. Crosses identify crashing
simulations (due to memory requirements becoming too severe).}
\centering
\resizebox{\textwidth}{!}{
\begin{tabular}{c|c|cc:cc:cc:cc:cc:cc:cc}
  Pb &$\tableIndices{\Delta x}{\Delta t}$& \multicolumn{2}{c:}{$2^{-1}$} & \multicolumn{2}{c:}{$2^{-2}$} & \multicolumn{2}{c:}{$2^{-3}$} & \multicolumn{2}{c:}{$2^{-4}$} 
                                         & \multicolumn{2}{c:}{$2^{-5}$} & \multicolumn{2}{c:}{$2^{-6}$} & \multicolumn{2}{c }{$2^{-7}$} \BotSp\\\hline
  \multirow{7}{*}{1}&$            2^{-2}$&$ 23$&$(23) $&$ 24$&$(24) $&$ 25$&$(25) $&$ 26$&$(26) $&$ 25$&$(25) $&$ 26$&$(26) $&$ 26$&$(26) $\TopSp\\
                    &$            2^{-3}$&$ 22$&$(22) $&$ 22$&$(22) $&$ 23$&$(23) $&$ 24$&$(26) $&$ 24$&$(26) $&$ 24$&$(24) $&$ 25$&$(25) $\\
                    &$            2^{-4}$&$ 22$&$(22) $&$ 22$&$(22) $&$ 23$&$(23) $&$ 23$&$(22) $&$ 22$&$(22) $&$ 23$&$(31) $&$ 22$&$(28) $\\
                    &$            2^{-5}$&$ 20$&$(20) $&$ 21$&$(21) $&$ 21$&$(21) $&$ 20$&$(20) $&$ 20$&$(20) $&$ 20$&$(28) $&$ 20$&$(21) $\\
                    &$            2^{-6}$&$ 19$&$(19) $&$ 19$&$(19) $&$ 19$&$(18) $&$ 19$&$(19) $&$ 19$&$(18) $&$ 19$&$(19) $&$ 20$&$(20) $\\
                    &$            2^{-7}$&$ 18$&$(18) $&$ 18$&$(18) $&$ 19$&$(18) $&$ 18$&$(18) $&$ 19$&$(19) $&$ 18$&$(18) $&$ 19$&$(19) $\\
                    &$            2^{-8}$&$ 17$&$(17) $&$ 18$&$(18) $&$ 18$&$(17) $&$ 17$&$(17) $&$ 18$&$(18) $&$ 17$&$(17) $&$ 16$&$(17) $\BotSp\\\hline
  \multirow{7}{*}{2}&$            2^{-2}$&$ 28$&$(28) $&$ 32$&$(32) $&$ 34$&$(34) $&$ 36$&$(36) $&$ 40$&$(40) $&$ 43$&$(43) $&$ 49$&$(50) $\TopSp\\
                    &$            2^{-3}$&$ 31$&$(31) $&$ 34$&$(34) $&$ 35$&$(35) $&$ 36$&$(41) $&$ 38$&$(42) $&$ 39$&$(38) $&$ 38$&$(39) $\\
                    &$            2^{-4}$&$ 30$&$(30) $&$ 32$&$(32) $&$ 33$&$(33) $&$ 34$&$(34) $&$ 34$&$(36) $&$ 35$&$(49) $&$ 35$&$(44) $\\
                    &$            2^{-5}$&$ 29$&$(29) $&$ 31$&$(30) $&$ 32$&$(32) $&$ 33$&$(33) $&$ 34$&$(35) $&$ 34$&$(49) $&$ 34$&$(37) $\\
                    &$            2^{-6}$&$ 28$&$(28) $&$ 30$&$(30) $&$ 31$&$(31) $&$ 32$&$(32) $&$ 32$&$(32) $&$ 33$&$(34) $&$ 31$&$(32) $\\
                    &$            2^{-7}$&$ 27$&$(27) $&$ 29$&$(28) $&$ 30$&$(30) $&$ 30$&$(30) $&$ 30$&$(30) $&$ 29$&$(30) $&$ 29$&$(30) $\\
                    &$            2^{-8}$&$ 25$&$(25) $&$ 26$&$(26) $&$ 27$&$(26) $&$ 28$&$(27) $&$ 26$&$(26) $&$ 26$&$(27) $&$ 26$&$(27) $\BotSp\\\hline
  \multirow{7}{*}{3}&$            2^{-2}$&$ 33$&$(34) $&$ 35$&$(40) $&$ 37$&$(42) $&$ 38$&$(38) $&$ 42$&$(42) $&$ 46$&$(46) $&$ 53$&$(55) $\TopSp\\
                    &$            2^{-3}$&$ 33$&$(33) $&$ 35$&$(35) $&$ 36$&$(37) $&$ 37$&$(49) $&$ 39$&$(47) $&$ 39$&$(39) $&$ 40$&$(40) $\\
                    &$            2^{-4}$&$ 31$&$(31) $&$ 33$&$(33) $&$ 34$&$(34) $&$ 35$&$(47) $&$ 35$&$(38) $&$ 36$&$(53) $&$ 37$&$(48) $\\
                    &$            2^{-5}$&$ 30$&$(30) $&$ 32$&$(32) $&$ 33$&$(33) $&$ 34$&$(34) $&$ 34$&$(35) $&$ 36$&$(53) $&$ 35$&$(39) $\\
                    &$            2^{-6}$&$ 28$&$(28) $&$ 31$&$(31) $&$ 32$&$(32) $&$ 32$&$(32) $&$ 33$&$(33) $&$ 33$&$(34) $&$ 33$&$(35) $\\
                    &$            2^{-7}$&$ 27$&$(27) $&$ 29$&$(29) $&$ 30$&$(30) $&$ 31$&$(31) $&$ 31$&$(31) $&$ 30$&$(30) $&$ 30$&$(30) $\\
                    &$            2^{-8}$&$ 25$&$(26) $&$ 27$&$(27) $&$ 28$&$(28) $&$ 28$&$(28) $&$ 28$&$(28) $&$ 26$&$(27) $&$\times$&$(\times)$\BotSp\\\hline
  \multirow{7}{*}{4}&$            2^{-2}$&$ 25$&$(25) $&$ 27$&$(27) $&$ 27$&$(27) $&$ 27$&$(27) $&$ 28$&$(28) $&$ 27$&$(27) $&$ 27$&$(27) $\TopSp\\
                    &$            2^{-3}$&$ 24$&$(24) $&$ 25$&$(25) $&$ 25$&$(25) $&$ 26$&$(28) $&$ 26$&$(27) $&$ 26$&$(26) $&$ 26$&$(26) $\\
                    &$            2^{-4}$&$ 24$&$(24) $&$ 25$&$(25) $&$ 25$&$(25) $&$ 25$&$(26) $&$ 24$&$(25) $&$ 24$&$(31) $&$ 23$&$(28) $\\
                    &$            2^{-5}$&$ 24$&$(24) $&$ 24$&$(26) $&$ 23$&$(25) $&$ 23$&$(24) $&$ 22$&$(24) $&$ 23$&$(32) $&$ 22$&$(24) $\\
                    &$            2^{-6}$&$ 21$&$(23) $&$ 22$&$(25) $&$ 22$&$(26) $&$ 21$&$(26) $&$ 21$&$(26) $&$ 22$&$(28) $&$ 21$&$(28) $\\
                    &$            2^{-7}$&$ 20$&$(24) $&$ 21$&$(26) $&$ 21$&$(28) $&$ 20$&$(28) $&$ 21$&$(31) $&$ 21$&$(32) $&$ 20$&$(33) $\\
                    &$            2^{-8}$&$ 19$&$(24) $&$ 20$&$(27) $&$ 20$&$(28) $&$ 19$&$(30) $&$ 20$&$(32) $&$ 19$&$(32) $&$ 19$&$(34) $\\
\end{tabular}}
\end{table}
In more detail, the application of the inverse of \cref{eqn::PT} involves
inverting three different operators: two of these are associated with the
application of the approximate space-time Schur complement \cref{eqn::XXdef},
and require inverting the pressure mass matrix $\mc{M}_p$
\cref{eqn::massStiffP::mass}, and the pressure ``Laplacian'' $\tilde{\mc{A}}_p$
\cref{eqn::Plaplacian}; the last one must deal with the monolithic space-time
matrix for the velocity variable $F_\bb{u}$ \cref{eqn::STFu}. In the ideal case,
we exploit LU solvers for both \cref{eqn::massStiffP::mass} and \cref{eqn::Plaplacian}, and
a sequential time-stepping procedure for \cref{eqn::STFu} (the latter further
requires inverting the spatial velocity operator $\mc{F}_{\bb{u},k}$
\cref{eqn::VRAD} at each time-step $k$, for which we also use LU).

However, resorting to direct methods becomes infeasible when solving systems of
large size, both in terms of memory and of computations required. Moreover, our
final goal is to actually parallelise in time the solution of the monolithic
system. In order to do so, we must substitute the time-stepping procedure for
the space-time velocity block $F_{\bb{u}}$ with a parallel alternative. Here, 
we use GMRES preconditioned by nonsymmetric AMG based on approximate ideal
restriction (AIR) \cite{lAIR,AIR}, as implemented in \textit{hypre}. AIR was
designed as a reduction-based classical AMG solver for advection-dominated
problems, but was shown to be effective on advection-diffusion problems as
well. More recently, AIR was demonstrated as an effective algebraic solver
for space-time finite element discretisations of advection(-diffusion)
\cite{sivas2020air}. Although likely not as effective as space-time geometric
multigrid techniques (e.g., \cite{Horton.1995}) for diffusive problems, AIR
is likely more robust for advective regimes and also less intrusive to apply, as
we simply need to form the space-time velocity matrix and pass it to
\textit{hypre}. Due to the use of an inner Krylov method, we swap the outer GMRES solver
with a \emph{flexible}-GMRES (FGMRES) solver \cite{FGMRES}. The inversions of the
pressure mass and stiffness matrices are also approximated: for $\mc{M}_p$, we
apply a fixed number of Chebyshev iterations \cite{andyCheb}, making use of the
bounds on its eigenvalues provided by \cite{andyEigBounds}; for
$\tilde{\mc{A}}_p$, we apply \textit{hypre}'s BoomerAMG \cite{AMG}.


\pgfplotsset{
  convResDiffTolPlotStyle/.style={
    tick label style={font=\small},
    xmin = 1e-13, xmax = 1e0, ymin = 1, ymax = 100,
    ymajorgrids= true,
    unbounded coords=jump,
    xlabel= Inner GMRES tolerance,
    ylabel= FGMRES iterations to convergence,
    ymajorgrids= true,
    axis background/.style={fill=gray!10},
    legend pos=north west,
    legend style={font=\small},
    cycle list = {{YlGnBu-D},{YlGnBu-F},{YlGnBu-H},{YlGnBu-L}},
  }
}
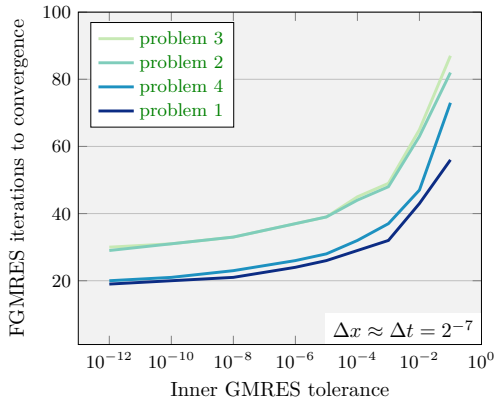
\begin{figure}[!t]
\centering
\resizebox{0.5\textwidth}{!}{
\begin{tikzpicture}[baseline,scale=0.8]
  \begin{semilogxaxis}[convResDiffTolPlotStyle]   
    \addplot+[line width=0.5mm] table[x index=0, y index=3] {data/convergence_results_Prec1_STsolve2_oU2_oP1_rc_SpaceTimeStokes_FGMRES_approx2_varFuiTol.txt};
    \addplot+[line width=0.5mm] table[x index=0, y index=2] {data/convergence_results_Prec1_STsolve2_oU2_oP1_rc_SpaceTimeStokes_FGMRES_approx2_varFuiTol.txt};
    \addplot+[line width=0.5mm] table[x index=0, y index=4] {data/convergence_results_Prec1_STsolve2_oU2_oP1_rc_SpaceTimeStokes_FGMRES_approx2_varFuiTol.txt};
    \addplot+[line width=0.5mm] table[x index=0, y index=1] {data/convergence_results_Prec1_STsolve2_oU2_oP1_rc_SpaceTimeStokes_FGMRES_approx2_varFuiTol.txt};
    \addlegendentry{\cref{pb::backStep}}
    \addlegendentry{\cref{pb::poiseuille}}
    \addlegendentry{\cref{pb::doubleGlazing}}
    \addlegendentry{\cref{pb::drivenCavity}}
    \node[anchor=south east, fill=white!10 ] at (axis cs: 1e0,1) {$\Delta x\approx\Delta t=2^{-7}$};
  \end{semilogxaxis}
\end{tikzpicture}}
\caption{\label{fig::FGMRESconvVarTol} Number of outer FGMRES iterations to convergence for different
tolerances of the inner GMRES solver for the space-time velocity block $F_\bb{u}$ in \cref{eqn::PT}.
The discretisation parameters are fixed to $\Delta x\approx\Delta t=2^{-7}$ for all problems in
\cref{sec::modelProblems}. The remaining options for the solvers involved are the same ones used
to fill the right column in \cref{tab::PconvIdealVSapprox}.}
\end{figure}

Results from using direct or approximate solvers are reported side-by-side
in \cref{tab::PconvIdealVSapprox}, for ease of comparison. The two cases show a
remarkably similar convergence behaviour, with the total number of iterations to
convergence remaining substantially unchanged. The only notable disparities are
observed for the most complex \cref{pb::doubleGlazing}, likely because
the 15 AIR iterations used as an approximate solver for the space-time velocity
block do not provide as accurate of an approximation as for \cref{pb::drivenCavity,pb::poiseuille,pb::backStep}.
Even here, the increase in iterations compared with exact inner solves is nicely
bounded, remaining below a factor of $2$ with respect to the ideal
case. This shows that the approximate
solvers employed can be as effective as their direct counterparts, and gives
an example of an effective time-parallel procedure which can work well in
tandem with the preconditioner proposed.

To draw a clearer picture of the sensitivity of our space-time block preconditioning
procedure to the approximate solution of the velocity block, we analyse more in
detail how the tolerance of the solver chosen for $F_{\bb{u}}$ impacts the
total number of iterations to convergence for the outer solver.
This gives some indication of the target accuracy required by the PinT solver for the
space-time block preconditioner to still be effective. From the results shown in
\cref{fig::FGMRESconvVarTol}, it transpires that we can afford solving the space-time
velocity system with a tolerance as lax as $\approx10^{-3}$, and still retain
convergence with a reasonably small number of iterations. Remarkably, this is in line
with the required tolerance for the velocity block solve $\mc{F}_{\bb{u},k}$ in
the single time-step preconditioner \cref{eqn::PTSingeStep} as observed,
for example, in \cite{southworth2020FP}.

We point out that, for the largest problems considered, the relevant systems
have sizes $N_u =526,338$, $N_p=66,049$, and $N_t=128$, for a total number of
unknowns of $(N_u+N_p)\cdot N_t\approx76\cdot10^6$; this further increases to
$N_u =5,775,362$, $N_p=722,945$, and $N_t=64$, totalling $\approx 416\cdot10^{6}$,
for \cref{pb::backStep} (a simulation with $\approx 832\cdot10^{6}$
degrees of freedom crashed due to lack of memory).
Compared to the noticeable size of the system, the
total number of iterations to convergence reported in
\cref{tab::PconvIdealVSapprox} remains small, regardless of mesh size and
problem type, providing evidence of the optimal scalability of the
preconditioners proposed. This is further confirmed by tracking more in detail
the evolution of the residual as the GMRES iterations progress: an example of
this is provided in \cref{fig::GMRESresEvol}, which shows reasonably uniform
convergence profiles, independently of the choice of both $\Delta x$ and $\Delta
t$.

\pgfplotsset{
  convResPlotStyle/.style={
    tick label style={font=\small},
    xmin = -1, xmax = 28, ymin = 1e-11, ymax = 1e1,
    ymajorgrids= true,
    unbounded coords=jump,
    xlabel= GMRES iteration,
    ylabel= Relative residual norm,
    ymajorgrids= true,
    axis background/.style={fill=gray!10},
    legend style={font=\small},
    cycle list = {{YlGnBu-D},{YlGnBu-E},{YlGnBu-F},{YlGnBu-G},{YlGnBu-H},{YlGnBu-J},{YlGnBu-L}},
  }
}
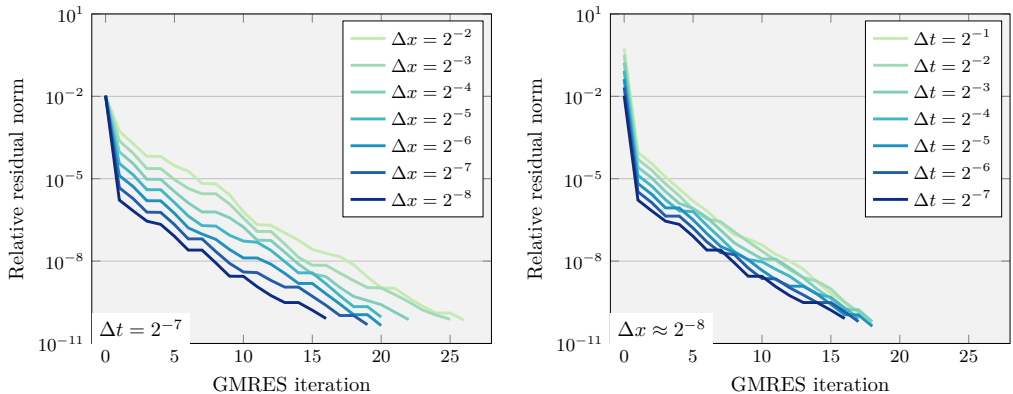
\begin{figure}[!t]
\begin{minipage}{0.49\textwidth}
  \centering
  \resizebox{\textwidth}{!}{
  \begin{tikzpicture}[baseline,scale=0.8]
    \begin{semilogyaxis}[convResPlotStyle]   
      \addplot+[line width=0.5mm] table[x index=0, y index=11]{data/Pb1_Prec1_STsolve0_oU2_oP1/NP128_r2.txt};
      \addplot+[line width=0.5mm] table[x index=0, y index=11]{data/Pb1_Prec1_STsolve0_oU2_oP1/NP128_r3.txt};
      \addplot+[line width=0.5mm] table[x index=0, y index=11]{data/Pb1_Prec1_STsolve0_oU2_oP1/NP128_r4.txt};
      \addplot+[line width=0.5mm] table[x index=0, y index=11]{data/Pb1_Prec1_STsolve0_oU2_oP1/NP128_r5.txt};
      \addplot+[line width=0.5mm] table[x index=0, y index=11]{data/Pb1_Prec1_STsolve0_oU2_oP1/NP128_r6.txt};
      \addplot+[line width=0.5mm] table[x index=0, y index=11]{data/Pb1_Prec1_STsolve0_oU2_oP1/NP128_r7.txt};
      \addplot+[line width=0.5mm] table[x index=0, y index=11]{data/Pb1_Prec1_STsolve0_oU2_oP1/NP128_r8.txt};
      \addlegendentry{$\Delta x=2^{-2}$}
      \addlegendentry{$\Delta x=2^{-3}$}
      \addlegendentry{$\Delta x=2^{-4}$}
      \addlegendentry{$\Delta x=2^{-5}$}
      \addlegendentry{$\Delta x=2^{-6}$}
      \addlegendentry{$\Delta x=2^{-7}$}
      \addlegendentry{$\Delta x=2^{-8}$}
      \node[anchor=south west, fill=white!10 ] at (axis cs: -1,1e-11) {$\Delta t=2^{-7}$};
    \end{semilogyaxis}
  \end{tikzpicture}}
\end{minipage}
\hfill
\begin{minipage}{0.49\textwidth}
  \centering
  \resizebox{\textwidth}{!}{
  \begin{tikzpicture}[baseline,scale=0.8]
    \begin{semilogyaxis}[convResPlotStyle]   
      \addplot+[line width=0.5mm] table[x index=0, y index=11]{data/Pb1_Prec1_STsolve0_oU2_oP1/NP2_r8.txt};
      \addplot+[line width=0.5mm] table[x index=0, y index=11]{data/Pb1_Prec1_STsolve0_oU2_oP1/NP4_r8.txt};
      \addplot+[line width=0.5mm] table[x index=0, y index=11]{data/Pb1_Prec1_STsolve0_oU2_oP1/NP8_r8.txt};
      \addplot+[line width=0.5mm] table[x index=0, y index=11]{data/Pb1_Prec1_STsolve0_oU2_oP1/NP16_r8.txt};
      \addplot+[line width=0.5mm] table[x index=0, y index=11]{data/Pb1_Prec1_STsolve0_oU2_oP1/NP32_r8.txt};
      \addplot+[line width=0.5mm] table[x index=0, y index=11]{data/Pb1_Prec1_STsolve0_oU2_oP1/NP64_r8.txt};
      \addplot+[line width=0.5mm] table[x index=0, y index=11]{data/Pb1_Prec1_STsolve0_oU2_oP1/NP128_r8.txt};
      \addlegendentry{$\Delta t=2^{-1}$}
      \addlegendentry{$\Delta t=2^{-2}$}
      \addlegendentry{$\Delta t=2^{-3}$}
      \addlegendentry{$\Delta t=2^{-4}$}
      \addlegendentry{$\Delta t=2^{-5}$}
      \addlegendentry{$\Delta t=2^{-6}$}
      \addlegendentry{$\Delta t=2^{-7}$}
      \node[anchor=south west, fill=white!10 ] at (axis cs: -1,1e-11) {$\Delta x\approx2^{-8}$};
    \end{semilogyaxis}
  \end{tikzpicture}}
\end{minipage}
\caption{\label{fig::GMRESresEvol} Details on relative (with respect to the norm
of the right-hand side) residual evolution of GMRES applied to
\cref{pb::drivenCavity}, using $P_T^{-1}$ with ideal components as a
right-preconditioner. Results on the left refer to simulations with fixed
$\Delta t=2^{-7}$ and varying $\Delta x$; results on the right are obtained by
fixing $\Delta x\approx2^{-8}$ and varying $\Delta t$. The set-up for the
solvers is similar to the one used for \cref{tab::PconvIdealVSapprox}.}
\end{figure}

\subsubsection{Dependence on P\'eclet number}
\label{sec::results::advection}
\begin{table}[b!]\footnotesize
\caption{\label{tab::PconvIdealPeclet} Number of iterations to convergence for
\cref{pb::doubleGlazing} with different values of \text{Pe}. Dashes represent
simulations that did not converge in the maximum prescribed number of iterations
($100$). Same solvers set-up as for the ideal case in \cref{tab::PconvIdealVSapprox}.}
\centering
\resizebox{\textwidth}{!}{
\begin{tabular}{c|ccccc:ccccc:ccccc:ccccc}
  $\tableIndices{\Delta x}{\Delta t}$&\multicolumn{5}{c:}{$2^{-4}$}&\multicolumn{5}{c:}{$2^{-5}$}
                                     &\multicolumn{5}{c:}{$2^{-6}$}&\multicolumn{5}{c }{$2^{-7}$}\\\hline
  Pe$\to $&$ 2^4$&$ 2^5$&$ 2^6$&$ 2^7$&$ 2^8 $&$  2^4$&$ 2^5$&$ 2^6$&$ 2^7$&$ 2^8 $&$  2^4$&$ 2^5$&$ 2^6$&$ 2^7$&$ 2^8 $&$  2^4$&$ 2^5$&$ 2^6$&$ 2^7$&$ 2^8 $\TopSp\\\hline
  $2^{-2}$&$ 32 $&$ 60 $&$ // $&$ // $&$ //  $&$  32 $&$ 57 $&$ // $&$ // $&$ //  $&$  31 $&$ 55 $&$ // $&$ // $&$ //  $&$  31 $&$ 54 $&$ // $&$ // $&$ //  $\TopSp\\
  $2^{-3}$&$ 28 $&$ 35 $&$ 54 $&$ // $&$ //  $&$  28 $&$ 35 $&$ 53 $&$ // $&$ //  $&$  28 $&$ 35 $&$ 53 $&$ // $&$ //  $&$  27 $&$ 35 $&$ 52 $&$ // $&$ //  $\\
  $2^{-4}$&$ 27 $&$ 32 $&$ 40 $&$ 55 $&$ //  $&$  27 $&$ 32 $&$ 39 $&$ 55 $&$ //  $&$  25 $&$ 31 $&$ 38 $&$ 53 $&$ //  $&$  25 $&$ 31 $&$ 38 $&$ 52 $&$ //  $\\
  $2^{-5}$&$ 25 $&$ 30 $&$ 35 $&$ 45 $&$ 66  $&$  25 $&$ 28 $&$ 35 $&$ 44 $&$ 66  $&$  24 $&$ 28 $&$ 34 $&$ 44 $&$ 64  $&$  23 $&$ 27 $&$ 34 $&$ 43 $&$ 63  $\\
  $2^{-6}$&$ 24 $&$ 27 $&$ 32 $&$ 39 $&$ 49  $&$  23 $&$ 27 $&$ 32 $&$ 38 $&$ 49  $&$  23 $&$ 27 $&$ 32 $&$ 37 $&$ 48  $&$  22 $&$ 26 $&$ 30 $&$ 37 $&$ 47  $\\
  $2^{-7}$&$ 23 $&$ 26 $&$ 31 $&$ 35 $&$ 40  $&$  22 $&$ 26 $&$ 30 $&$ 35 $&$ 41  $&$  22 $&$ 26 $&$ 29 $&$ 33 $&$ 41  $&$  21 $&$ 25 $&$ 29 $&$ 33 $&$ 40  $\\
  $2^{-8}$&$ 22 $&$ 25 $&$ 29 $&$ 32 $&$ 37  $&$  21 $&$ 25 $&$ 29 $&$ 32 $&$ 37  $&$  21 $&$ 25 $&$ 28 $&$ 32 $&$ 37  $&$  20 $&$ 24 $&$ 27 $&$ 32 $&$ 37  $
\end{tabular}}
\end{table}
When applied to the solution of Oseen equations, we observe that the performance
of the preconditioner degrades as we increase the intensity of the advection field
imposed, as foreshadowed by the analysis in \cref{sec::precon::eigs}.
This is shown in \cref{tab::PconvIdealPeclet}, where we progressively ramp up the value 
of Pe in \cref{pb::doubleGlazing}, solving for various $\Delta x$: convergence fails
for coarser spatial meshes and largest P\'eclet numbers. This hints at the necessity
of either employing appropriately refined meshes, or relying on stabilisation techniques
\cite[Chap.~13.8]{quarterNM} (which we did not apply in our analysis), whenever
advection-dominated flows need to be resolved. The preconditioner for the single time-step
case \cref{eqn::PTSingeStep} was shown to suffer from a similar lack of robustness
\cite[Table~9.3]{andyFIT}, so it is somewhat expected that this weakness carries over to the
whole space-time case.
However, we point out that the performance of the preconditioner does not vary significantly as long as
we keep constant the \emph{grid} P\'eclet number $\tilde{\text{Pe}} \coloneqq \Delta x \text{Pe} / L$
\cite[Chap.~13.2]{quarterNM}. Unlike its \emph{global} counterpart \cref{eqn::peclet},
this parameter provides a \emph{local} measure of the dominance of advection over diffusion,
and is useful in identifying conditions for which instabilities and oscillations might
arise in the numerical solution. This parameter is kept fixed down each diagonal of 
\cref{tab::PconvIdealPeclet}, which show only slow growth in the number of iterations
to convergence with increasing Pe.

\subsection{The nonlinear setting: incompressible Navier-Stokes}
\label{sec::results::NS}
As stated in \cref{sec::oseen}, Oseen equations \cref{eqn::oseen} can be interpreted
as a linearisation of the Navier-Stokes system \cref{eqn::navierStokes}.
If we then introduce an \emph{outer} solver to tackle the nonlinear advection
term in Navier-Stokes, we can use our preconditioner to accelerate the solution
of the corresponding \emph{inner} (linearised) problems, much like what is done
for the single time-step case in \cite[Chap.~10]{andyFIT}. We experiment on the
performance of our preconditioner $P_T$ \cref{eqn::PT} in this situation, and
resolve the nonlinearity in
\cref{eqn::navierStokes} via \emph{Picard iterations} \cite{picard}. That is, we
iteratively solve the linearised system 
\begin{equation}
  \left[\begin{array}{c|c}
    F_{\bb{u}}((\bb{u})^{j}) & B^T\BotSp\\\hline
    B 								       &   \mathbf{0} \TopSp	
  \end{array}\right]\left[\begin{array}{c}
    (\bb{u})^{j+1}\BotSp\\\hline
    (\bb{p})^{j+1}\TopSp
  \end{array}\right]=\left[\begin{array}{c}
    \bb{f}\BotSp\\\hline
    \bb{0}\TopSp
  \end{array}\right],
  \label{eqn::picard}
\end{equation}
until convergence is reached. Here, $(\bb{u})^{j+1}$ and $(\bb{p})^{j+1}$ denote
the whole space-time velocity and pressure solutions at the $j+1$-th Picard iteration,
while we made explicit the dependence of the space-time velocity block
$F_{\bb{u}}=F_{\bb{u}}(\bb{w})$ on the advection field, which we substitute with
the approximate solution at the previous iteration $(\bb{u})^j$. 

\begin{table}[b!]\footnotesize
\caption{\label{tab::PconvIdealNS} Convergence results for Navier-Stokes with the setup
of \cref{pb::drivenCavity,pb::backStep}. Picard iterations used as
the outer nonlinear solver, which achieves convergence with a tolerance of
$10^{-9}$. Inside each column, on the left is reported the number of outer
iterations to convergence, while on the right (in brackets) the average number
of inner GMRES iterations per outer Picard iteration. The inner GMRES solver has
the options used for recovering the results in the left column of
\cref{tab::PconvIdealVSapprox}. Notice \cref{pb::poiseuille} is missing: this is
because the nonlinearity does not affect the analytical solution for
Poiseuille, and hence results are not particularly informative in this case
(convergence occurs at the very first Picard iteration in almost the totality of
the experiments). Also \cref{pb::doubleGlazing} is missing, since it already
represents a linearisation of Navier-Stokes for
\cref{pb::drivenCavity}. Crosses identify crashing
simulations (due to memory requirements becoming too severe).}
\centering
\resizebox{\textwidth}{!}{
\begin{tabular}{c|c|cc:cc:cc:cc:cc:cc:cc}
  Pb&$\tableIndices{\Delta x}{\Delta t}$& \multicolumn{2}{c:}{$2^{-1}$} & \multicolumn{2}{c:}{$2^{-2}$} & \multicolumn{2}{c:}{$2^{-3}$} & \multicolumn{2}{c:}{$2^{-4}$} 
                                        & \multicolumn{2}{c:}{$2^{-5}$} & \multicolumn{2}{c:}{$2^{-6}$} & \multicolumn{2}{c }{$2^{-7}$} \BotSp\\\hline
  \multirow{7}{*}{1}&$ 2^{-2}$&$ 5$&$(11.60)  $&$ 4$&$(14.75) $&$ 4$&$(14.25) $&$ 4$&$(14.50) $&$ 4$&$(13.50) $&$ 4$&$(13.25) $&$ 4$&$(13.00)  $\TopSp\\
                    &$ 2^{-3}$&$ 4$&$(13.50)  $&$ 4$&$(13.25) $&$ 4$&$(13.25) $&$ 4$&$(13.00) $&$ 4$&$(12.50) $&$ 4$&$(12.25) $&$ 5$&$( 9.40)  $\\
                    &$ 2^{-4}$&$ 4$&$(12.50)  $&$ 4$&$(12.50) $&$ 4$&$(12.25) $&$ 4$&$(12.00) $&$ 4$&$(11.50) $&$ 5$&$( 8.80) $&$ 5$&$( 8.60)  $\\
                    &$ 2^{-5}$&$ 4$&$(11.50)  $&$ 4$&$(11.00) $&$ 4$&$(11.50) $&$ 4$&$(10.75) $&$ 5$&$( 8.00) $&$ 5$&$( 8.20) $&$ 4$&$( 9.25)  $\\
                    &$ 2^{-6}$&$ 4$&$(10.50)  $&$ 4$&$(10.00) $&$ 4$&$(10.00) $&$ 4$&$( 9.50) $&$ 4$&$( 9.25) $&$ 4$&$( 8.50) $&$ 4$&$( 8.75)  $\\
                    &$ 2^{-7}$&$ 4$&$( 9.75)  $&$ 4$&$( 9.00) $&$ 4$&$( 9.25) $&$ 5$&$( 6.80) $&$ 4$&$( 8.25) $&$ 4$&$( 8.00) $&$ 4$&$( 7.75)  $\\
                    &$ 2^{-8}$&$ 3$&$(11.00)  $&$ 3$&$(11.00) $&$ 3$&$(10.67) $&$ 3$&$(10.33) $&$ 3$&$(10.33) $&$ 4$&$( 7.25) $&$ 4$&$( 6.75)  $\BotSp\\\hline
  \multirow{7}{*}{3}&$ 2^{-2}$&$ 5$&$(19.80) $&$ 5$&$(18.20) $&$ 5$&$(18.00) $&$ 5$&$(18.00) $&$ 5$&$(18.40) $&$ 5$&$(19.00) $&$ 5$&$(21.40) $\TopSp\\
                    &$ 2^{-3}$&$ 5$&$(21.80) $&$ 5$&$(19.60) $&$ 5$&$(17.80) $&$ 5$&$(17.00) $&$ 5$&$(16.80) $&$ 5$&$(16.40) $&$ 5$&$(16.20) $\\
                    &$ 2^{-4}$&$ 5$&$(22.60) $&$ 5$&$(20.60) $&$ 5$&$(18.80) $&$ 5$&$(16.80) $&$ 5$&$(15.60) $&$ 5$&$(14.80) $&$ 5$&$(14.20) $\\
                    &$ 2^{-5}$&$ 5$&$(26.60) $&$ 5$&$(22.80) $&$ 5$&$(20.60) $&$ 5$&$(17.60) $&$ 5$&$(16.20) $&$ 5$&$(15.80) $&$ 5$&$(14.20) $\\
                    &$ 2^{-6}$&$ 4$&$(36.00) $&$ 4$&$(31.75) $&$ 4$&$(28.25) $&$ 5$&$(18.00) $&$ 5$&$(17.20) $&$ 5$&$(16.00) $&$ 5$&$(14.80) $\\
                    &$ 2^{-7}$&$ 4$&$(43.00) $&$ 4$&$(37.00) $&$ 4$&$(31.75) $&$ 4$&$(24.75) $&$ 5$&$(18.40) $&$ 5$&$(18.60) $&$ 5$&$(15.00) $\\
                    &$ 2^{-8}$&$ 4$&$(52.50) $&$ 4$&$(46.25) $&$ 4$&$(33.50) $&$ 4$&$(28.25) $&$ 5$&$(22.40) $&$\times$&$(\times)$&$\times$&$(\times)$
\end{tabular}}
\end{table}

Linear and nonlinear convergence obtained by applying a Picard linearisation
\cref{eqn::picard} to the problems in \cref{sec::modelProblems} are reported in
\cref{tab::PconvIdealNS}. Especially for \cref{pb::drivenCavity}, both the inner
and outer solve iterations remain $\mathcal{O}(1)$, independently on the spatial
and temporal mesh spacing. Similar results holds for the inflow-outflow problem
\cref{pb::backStep}, but here the average number of inner GMRES iterations is
larger. In the specific case of a very large time-step relative to the spatial
mesh, the number of inner iterations can be a fair amount larger. However, this
is not necessarily surprising, as \cref{pb::backStep} offers more complex dynamical
behavior and we are resolving the nonlinearity over a very large time-step. This
likely results in a very stiff and ill-conditioned linear system to solve for
each nonlinear iteration. Nevertheless, the outer space-time Picard linearisation
remains perfectly scalable.

\subsection{Comparison with sequential time-stepping}
\label{sec::results::STvsSTS}
\Cref{sec::precon::STvsTS} explains how the two
main factors in analysing the effectiveness of our preconditioning strategy
compared with classical time-stepping are given by the ratios
$C_{F_{\bb{u}}}/(C_{\mc{F}_{\bb{u}}}\,N_t)$ and $N_{it}^{ST}/N_{it}^0$. The
first ratio directly refers to the efficiency of the method chosen for the
parallel-in-time integration of the velocity block $F_{\bb{u}}$: identifying the
optimal solver in this sense is still a matter of research, and a detailed
investigation comparing the effectiveness of the various parallel-in-time
methods available in the literature \cite{PinTgander} goes well beyond the scope
of this paper. Here, we focus on the second ratio because it is
intrinsic to the preconditioning strategy proposed: in fact, it describes the
overhead produced by tackling the whole space-time system at once, rather than
via sequential time-stepping.

\begin{table}[b!]\footnotesize
\caption{\label{tab::PCompSingleStep} Ratio $N_{it}^{ST}/N_{it}^0$, where
$N_{it}^{ST}$ is taken from the ideal case in \cref{tab::PconvIdealVSapprox},
while $N_{it}^0$ represents the average number of iterations to convergence per
time-step, solving for \cref{eqn::STStepDisc} via forward substitution. For the
latter, GMRES is used at each time-step, right-preconditioned with the single
time-step counterpart of $P_T$, \cref{eqn::PTSingeStep} with
\cref{eqn::XXdefSingeStep}. Convergence is reached with a tolerance of
$10^{-10}/\sqrt[]{N_t}$. The remaining options for the solvers involved are the
same ones used to fill the left column in \cref{tab::PconvIdealVSapprox}. Crosses identify crashing
simulations (due to memory requirements becoming too severe).}
\centering
\resizebox{\textwidth}{!}{
\begin{tabular}{c|c|cccccc||cccccc|c}
  Pb &$\tableIndices{\Delta x}{\Delta t}$&$2^{-2}$&$2^{-3}$&$2^{-4}$&$2^{-5}$&$2^{-6}$&$2^{-7}$
  																			 &$2^{-2}$&$2^{-3}$&$2^{-4}$&$2^{-5}$&$2^{-6}$&$2^{-7}$&Pb\BotSp\\\hline
  \multirow{7}{*}{1}&$            2^{-2}$&$1.17$&$1.32$&$1.33$&$1.39$&$1.58$&$1.78  $&$  1.39$&$1.55$&$1.76$&$2.11$&$2.43$&$2.98$&\multirow{7}{*}{3}\TopSp\\
                    &$            2^{-3}$&$1.02$&$1.08$&$1.23$&$1.45$&$1.60$&$1.74  $&$  1.28$&$1.42$&$1.64$&$1.81$&$2.00$&$2.27$&                        \\
                    &$            2^{-4}$&$1.02$&$1.16$&$1.30$&$1.33$&$1.46$&$1.51  $&$  1.21$&$1.31$&$1.44$&$1.58$&$1.80$&$2.09$&                        \\
                    &$            2^{-5}$&$1.06$&$1.17$&$1.15$&$1.21$&$1.27$&$1.34  $&$  1.17$&$1.28$&$1.42$&$1.56$&$1.82$&$2.07$&                        \\
                    &$            2^{-6}$&$1.03$&$1.07$&$1.09$&$1.14$&$1.21$&$1.33  $&$  1.18$&$1.27$&$1.37$&$1.57$&$1.86$&$2.28$&                        \\
                    &$            2^{-7}$&$1.00$&$1.09$&$1.04$&$1.17$&$1.16$&$1.33  $&$  1.15$&$1.26$&$1.45$&$1.70$&$1.98$&$2.35$&                        \\
                    &$            2^{-8}$&$1.03$&$1.05$&$1.03$&$1.16$&$1.20$&$1.27  $&$  1.16$&$1.29$&$1.52$&$1.74$&$1.92$&$\times$&                \BotSp\\\hline
  \multirow{7}{*}{2}&$            2^{-2}$&$1.42$&$1.62$&$1.83$&$2.12$&$2.53$&$3.13  $&$  1.24$&$1.29$&$1.34$&$1.48$&$1.59$&$1.81$&\multirow{7}{*}{4}\TopSp\\
                    &$            2^{-3}$&$1.36$&$1.44$&$1.61$&$1.82$&$2.08$&$2.26  $&$  1.10$&$1.12$&$1.25$&$1.46$&$1.61$&$1.73$&                        \\
                    &$            2^{-4}$&$1.19$&$1.26$&$1.34$&$1.49$&$1.74$&$1.96  $&$  1.12$&$1.16$&$1.30$&$1.33$&$1.42$&$1.41$&                        \\
                    &$            2^{-5}$&$1.15$&$1.20$&$1.28$&$1.48$&$1.66$&$1.99  $&$  1.19$&$1.17$&$1.21$&$1.22$&$1.30$&$1.29$&                        \\
                    &$            2^{-6}$&$1.10$&$1.17$&$1.28$&$1.44$&$1.81$&$2.19  $&$  1.10$&$1.10$&$1.11$&$1.17$&$1.22$&$1.29$&                        \\
                    &$            2^{-7}$&$1.10$&$1.20$&$1.30$&$1.57$&$1.89$&$2.46  $&$  1.06$&$1.05$&$1.06$&$1.19$&$1.25$&$1.26$&                        \\
                    &$            2^{-8}$&$1.07$&$1.19$&$1.42$&$1.58$&$1.95$&$2.54  $&$  1.05$&$1.05$&$1.05$&$1.19$&$1.21$&$1.37$&                  \BotSp\\
\end{tabular}}
\end{table}

To give an indication of how $N_{it}^0$ and $N_{it}^{ST}$ compare with each
other when ideal solvers are employed, we collect convergence results from the
application of \cref{eqn::PTSingeStep} within the time-stepping procedure, for
problems equivalent to those analysed in \cref{tab::PconvIdealVSapprox}. To
render the comparison more fair, for the single time-step we ask for a stricter
tolerance on convergence, scaled by a factor $\sqrt[]{N_t}$ with respect to the
whole space-time case: in this way, the Euclidean norm of the final space-time
residual is similar in the two cases. We report the resulting ratio
$N_{it}^{ST}/N_{it}^0$ in \cref{tab::PCompSingleStep}.
We can see that this sits between
roughly $1$ and $3$, in most cases closer to one, but becoming larger as $\Delta t\to0$,
particularly $\Delta t \ll \Delta x$. This observation can be
explained by pointing out that most of the advantage of time-stepping is given
by the fact that, at each iteration, a good initial guess is provided by the
value of the solution at the previous time-step, which we exploit in our
experiments: with a reduced step-size, it is safe to expect that also the
solution varies less between one time-step and the next. Nevertheless,
results demonstrate little-to-no degradation in convergence speed
between the two approaches and, thus, minimal overhead in applying block
preconditioning techniques in an all-at-once fashion, compared with
sequential time-stepping.

\end{section}
\begin{section}{Conclusion and future work}
\label{sec::conclusion}
In this work, we introduce a novel all-at-once block preconditioner for the
discretised time-dependent Oseen equations. We take inspiration
from heuristics that have proven successful in the single time-step or steady-state
 settings, and apply similar principles to the space-time setting. The preconditioning procedure is designed to
better expose the system to effective strategies for time-parallelisation by
partitioning the preconditioning into two stages: (i) the space-time solution of
a time-dependent advection-diffusion equation, and (ii) a (spatial) mass-matrix
inverse and pressure Laplacian inverse at each time point. The latter can easily
be applied in a fast parallel manner using existing methods and is nearly
analogous to sequential time-stepping. Stage (i) reduces the solution of a
space-time system of PDEs to the (approximate) solution of a single-variable
space-time PDE, a much more tractable problem for fast, parallel
solvers. In particular, a number of approaches have already
seen success here, including space-time geometric multigrid
\cite{Falgout.2017,Horton.1995} and AMG applied to the space-time matrix \cite{sivas2020air},
as well as the more standard PinT methods \cite{PinTgander}.

The approach pursued here provides a different take on time-parallelisation for
time-dependent incompressible flows, compared with what is typically done in the
literature \cite{miao2019convergence,trindade2006parallel,trindade2004parallel,
croce2014parallel,pararealNS}. To our knowledge, most of the latter leaves the
underlying time-stepping procedure for solving \cref{eqn::STStepDisc}
largely unchanged, and rather focuses on building an overarching
framework which breaks its innate sequentiality. Instead, our method exploits
the block structure of the system \emph{at the space-time level}, in order to
simplify the type of problem we need to (parallel-)integrate in time. We can
draw an analogy between the role of the preconditioner for the single time-step
case \cref{eqn::PTSingeStep} and our space-time preconditioner \cref{eqn::PT}:
the first provides an effective method for solving Oseen via sequential time
stepping, so long as a fast solver for an elliptic operator is available; the
second provides an effective method for solving \textit{space-time} Oseen, so
long as a fast (time-parallel) solver for a \emph{parabolic} operator is
available.

The effectiveness of our preconditioner is demonstrated by applying it to a
variety of model problems taken from the literature, which cover a range of flow
configurations. The main measure for performance is given by the number of
iterations required for convergence, which we have shown to be scalable in
spatial and temporal mesh size, and comparable to those required for a classical
time-stepping (non-time-parallel) routine. A simple extension to a nonlinear
Navier-Stokes problem is investigated as well, which is linearised via a Picard
iteration, wherein we observe near perfect scalability of linear and nonlinear
iterations.

Future research will revolve around developing fast, robust, and time-parallel
solution methods for the scalar space-time velocity block in our preconditioner,
and extending the principles developed in this paper to the space-time solution
of other time-dependent systems of PDEs.


\section*{Acknowledgments}

Los Alamos National Laboratory report number LA-UR-21-20115.

\end{section}

\bibliographystyle{siamplain}
\bibliography{ex_article.bib}
\end{document}